\documentclass[12pt]{article}%
\usepackage{amsmath}
\usepackage{mitpress}
\usepackage{graphicx}%

\usepackage{epsfig}

\usepackage{amsfonts}%
\usepackage{amssymb}
\newtheorem{theorem}{Theorem}

\newtheorem{conjecture}[theorem]{Conjecture}
\newtheorem{corollary}[theorem]{Corollary}

\newtheorem{definition}[theorem]{Definition}
\newtheorem{example}[theorem]{Example}

\newtheorem{lemma}[theorem]{Lemma}

\newtheorem{problem}[theorem]{Problem}
\newtheorem{proposition}[theorem]{Proposition}
\newtheorem{remark}[theorem]{Remark}

\newenvironment{proof}[1][Proof]{\textbf{#1.} }{\ \rule{0.5em}{0.5em}}
\newdimen\dummy
\dummy=\oddsidemargin
\begin{document}

\title{Metric Coordinate Systems}
\author{{\small Craig Calcaterra, Axel Boldt, Michael Green}\\{\small Department of Mathematics; Metropolitan State University;}\\{\small St. Paul, MN 55106;\bigskip\ craig.calcaterra@metrostate.edu}\\{\small David Bleecker}\\{\small Department of Mathematics; University of Hawaii;}\\{\small Honolulu, HI 96822}}
\maketitle

\begin{abstract}
Coordinate systems are defined on general metric spaces with the purpose of
generalizing vector fields on a manifold. Conversion formulae are available
between metric and Cartesian coordinates on a Hilbert space. Nagumo's
Invariance Theorem is invoked to prove the analogue of the classical
Cauchy-Lipschitz Theorem for vector fields on a locally compact coordinatized
space. A metric space version of Nagumo's Theorem is one consequence. Examples
are given throughout.



\end{abstract}

\section{Introduction}

The notion of a metric coordinate system is offered here to extend the methods
of calculus and differential equations to metric spaces. The inspiration
behind metric coordinates is quite simple. On the plane $\mathbb{E}^{2}$, for
example, choose three non-colinear points $a,$ $b,$ and $c$. Then every point
$x\in\mathbb{E}^{2}$ is distinguished by three numbers, $d\left(  x,a\right)
,$ $d\left(  x,b\right)  ,$ and $d\left(  x,c\right)  $, which we call the
metric coordinates of $x.$ In a similar manner any metric space may be coordinatized.

The idea of metric coordinates has been put forth in the past to study static
problems in Euclidean spaces: \cite{Ellis}, \cite{Goldberg}, \cite{Havel}.
There have been several recent and notable efforts to develop generalizations
of differential equations in the context of metric spaces: quasi-differential
equations \cite{Panasyuk1}, mutational analysis \cite{Aubin}, and arc fields
\cite{CalcNBleecker}. These largely commensurable approaches have each
succeeded in producing a generalization of the Cauchy-Lipschitz Theorem. In
each of these schemes the idea of velocity at a point $x$ in a metric space
$\left(  X,d\right)  $ is represented by a curve issuing from $x$. The method
of this paper is different.

With metric coordinate systems, $X$ is embedded into a Banach subspace $E$
of$\ \mathbb{R}^{C}$ where $C$ is the set of coordinatizing points. $E$ is
used to define vector fields on $X$. Under suitable assumptions, the vector
field can be extended to a Lipschitz continuous vector field on $E$. Then the
traditional Cauchy-Lipschitz Theorem on Banach spaces yields unique solutions.
The Nagumo Invariance Theorem then promises that solutions with initial
conditions in the embedded subset remain there. The proof allows arbitrary
coordinatizing sets $C$, but uses local compactness. We expect a version
without this restriction is possible.

One of the strengths of metric coordinate systems is that, due to the
embedding, solving metric-coordinate vector fields on $X$ reduces to solving
an ODE on $\mathbb{R}^{C}$. Our vector fields and solutions will depend on the
choice of $C$. This coordinate dependence may be an advantage because it
allows us to capture dynamics that cannot easily be described otherwise. Also
metric coordinates, like other types of coordinate systems, are often more
convenient than Cartesian coordinates for solving certain problems. Spheres,
ellipses and hyperbolae are the loci of linear equations in metric coordinates.

Throughout the paper, examples are explored on Euclidean and non-Euclidean
spaces. Several open lines of research are detailed in the concluding section.

\section{Metric coordinatizing sets}

\begin{definition}
Let $\left(  M,d\right)  $ be a metric space with $X\subset M$. A
\textbf{metric coordinatizing set} for $X$ is a set of points $C\subset M$
with the property that for all $x,y\in X$ with $x\neq y$, there is some $c\in
C$ such that $d\left(  x,c\right)  \neq d\left(  y,c\right)  .$

We then call $\left(  M,d,X,C\right)  $ a \textbf{metric coordinate system}.
\end{definition}

As any point $x\in X$ in a metric coordinate system $\left(  M,d,X,C\right)  $
is represented by a $C$-tuple of real numbers $x_{C}=\left(  x_{c}\right)
_{c\in C}$, this will be called the $C$\textbf{\ embedding} of $X$ into
$\mathbb{R}^{C}.$ We are using the term ``embedding'' loosely here; it is not
necessarily a homeomorphism onto its image as the inverse is not necessarily continuous.

Throughout the paper we will be using arbitrary sets $C\subset M$ which may be
infinite or even unbounded. Most examples, however, suffice with finite sets
$C$ as in the following:

\begin{example}
\label{H^2CoordEx}We begin with Euclidean spaces. Consider the open half-plane
$H^{2}$ in the Euclidean plane $\mathbb{E}^{2}$ with the Euclidean metric $d$.
Pick any two distinct points $a$ and $b$ on the boundary. We can locate any
point $x$ in $H^{2}$ if we know its distances to $a$ and $b$, say
$x_{a}=d\left(  x,a\right)  $ and $x_{b}=d\left(  x,b\right)  $. Then $\left(
\mathbb{E}^{2},d,H^{2},\left\{  a,b\right\}  \right)  $ is a
\textit{bona-fide} metric coordinate system.

Equations in $\left(  \mathbb{E}^{2},d,H^{2},\left\{  a,b\right\}  \right)  $
are naturally different from those in Cartesian or polar coordinate systems.
E.g., for any $r>d\left(  a,b\right)  ,$ the locus of the equation%
\begin{equation}
x_{a}+x_{b}=r \label{x+y=c}%
\end{equation}
in metric coordinates is the set
\[
\left\{  x\in H^{2}:d\left(  x,a\right)  +d\left(  x,b\right)  =r\right\}  .
\]
The graph of $\left(  \ref{x+y=c}\right)  $ is half of an ellipse with foci at
$a$ and $b.$

$\mathbb{E}^{2}$, the plane, requires 3 non-colinear points for a metric
coordinatizing set. $H^{3}$ (the half-space) is metrically coordinatized with
3 non-colinear points on its boundary, and $\mathbb{E}^{3}$ needs 4
non-coplanar points. Many geometrical objects are readily described in metric
coordinates on $\mathbb{E}^{3}$:%
\[%
\begin{array}
[c]{ll}%
\text{Sphere (center }a,\text{ radius }r\text{)} & x_{a}=r\qquad
r\geq0\medskip\\
\text{Ellipsoid (foci }a,b\text{)} & x_{a}+x_{b}=r\text{\qquad}r\geq d\left(
a,b\right)  \medskip\\
\text{Hyperboloid (foci }a,b\text{)} & \left|  x_{a}-x_{b}\right|
=r\text{\qquad}0<r<d\left(  a,b\right)  \medskip\\
\text{Infinite Cylinder} & \sqrt{s\left(  s-x_{a}\right)  \left(
s-x_{b}\right)  \left(  s-d\left(  a,b\right)  \right)  }=r\medskip\\
\text{ (with axis }\overleftrightarrow{ab}\text{, radius }\frac{2r}{d\left(
a,b\right)  }\text{)} & \qquad\text{where }s=\frac{x_{a}+x_{b}+d\left(
a,b\right)  }{2}\medskip\\
\text{Infinite Cone} & x_{b}^{2}=d\left(  a,b\right)  ^{2}+x_{a}^{2}%
-2x_{a}d\left(  a,b\right)  \cos\theta\medskip\\
\text{ (with axis }\overleftrightarrow{ab}\text{, vertex }a\text{, angle
}\theta\text{)} & \medskip\\
\text{Plane (}\bot\overleftrightarrow{ab}\text{)} & x_{a}=x_{b}\medskip\\
\text{Segment }\overline{ab} & x_{a}+x_{b}=d\left(  a,b\right)  \medskip\\
\text{Ray }\overrightarrow{ab} & x_{a}\pm x_{b}=d\left(  a,b\right)
\medskip\\
\text{Line }\overleftrightarrow{ab} & \left|  x_{a}\pm x_{b}\right|  =d\left(
a,b\right)  \medskip
\end{array}
\]
The equation for the cylinder comes from Heron's formula for area of a
triangle. The equation for the cone is simply the cosine angle formula for a
triangle and represents only one half of a two sided cone; the other half is
given when $\theta$ is replaced with $\pi-\theta$. More general equations for
lines and planes are available but are not so concise. Choosing the
coordinates according to the problem simplifies the formulae.

As each of the above formulae use only metric coordinates, they may serve as
definitions for the various geometrical objects in general metric spaces.
\end{example}

\begin{proposition}
\label{metTriangleProp}Let $\left(  M,d,X,C\right)  $ be a metric coordinate
system. The metric coordinates $x_{c}:=d\left(  x,c\right)  $ of any point
$x\in X$ satisfy%
\begin{align}
x_{c}  &  \geq0\nonumber\\
\left|  x_{a}-x_{b}\right|   &  \leq d\left(  a,b\right) \label{MetTriangle}\\
x_{a}+x_{b}  &  \geq d\left(  a,b\right) \nonumber
\end{align}
for all $a,b,c\in C$.
\end{proposition}

\begin{proof}
Triangle inequality.
\end{proof}

Proposition \ref{metTriangleProp}, though mathematically trivial is used in
every example. It shows that care must be taken when defining curves in terms
of metric coordinates since not all $C$-tuples describe points in $X$.

\begin{example}
On any metric space $\left(  M,d\right)  $ there is at least one metric
coordinatizing set for any subset $X$. The worst-case scenario is the discrete
metric, defined on any set $M$ as%
\[
d\left(  x,y\right)  :=\left\{
\begin{array}
[c]{cc}%
1 & \text{ if }x\neq y\\
0 & \text{ if }x=y.
\end{array}
\right.
\]
This metric requires all of the points in $X$ save one for its metric
coordinatizing set.
\end{example}

\begin{example}
On a separable metric space $M$, any subset $X$ may be coordinatized with
countably many points.
\end{example}

\begin{example}
\label{InfinityMetric}Take $M=X=\mathbb{R}^{2}$ with the supremum metric
\[
d_{\infty}\left(  x,y\right)  :=\underset{i=1,2}{\max}\left\{  \left|
x_{i}-y_{i}\right|  \right\}  \text{.}%
\]
No bounded set is a coordinatizing set for $X$. If $C$ is contained in some
square, then two vertically aligned points placed far enough to the left of
the square cannot be distinguished by $C$. $X$ may, however, be coordinatized
by%
\[
C:=\left(  \mathbb{N\times}\left\{  0\right\}  \right)  \cup\left(  \left\{
0\right\}  \mathbb{\times N}\right)  .
\]
where $\mathbb{N}:=\left\{  1,2,3,...\right\}  $
\end{example}

\section{Conversion formulae for Hilbert spaces}

On the Euclidean plane $\mathbb{E}^{2}$ choose metric coordinates $a,b,c$ so
the rays $\overrightarrow{ca}$ and $\overrightarrow{cb}$ are perpendicular
with $d\left(  a,c\right)  =1=d\left(  b,c\right)  $. Define a Cartesian
coordinate system on the plane with the origin $\left(  0,0\right)  $ at $c,$
the positive $x$-axis along the ray $\overrightarrow{ca}$ and the positive $y
$-axis along the ray $\overrightarrow{cb}.$ The conversion
formulae$\footnote{To write $\left(  w_{a},w_{b},w_{c}\right)  =w=\left(
x,y\right)  $ is technically abuse of notation. $\left(  w_{a},w_{b}%
,w_{c}\right)  $ and $\left(  x,y\right)  $ are actually representations of
$w,$ and in the sequel we write $w_{C}=\left(  w_{a},w_{b},w_{c}\right)  $ to
make this distinction explicit.}$ are easy to find:%
\begin{gather}
\text{Metric}\quad\quad\left(  w_{a},w_{b},w_{c}\right)  =w=\left(
x,y\right)  \quad\quad\text{Cartesian}\label{3metCoordEq}\\
w_{c}=\sqrt{x^{2}+y^{2}}\nonumber\\
w_{b}=\sqrt{x^{2}+\left(  y-1\right)  ^{2}}\nonumber\\
w_{a}=\sqrt{\left(  x-1\right)  ^{2}+y^{2}}.\nonumber
\end{gather}
Solving these same equations for $x$ and $y$ yields the inverse formulae%
\begin{equation}
x=\dfrac{w_{c}^{2}-w_{a}^{2}+1}{2}\text{\qquad and\qquad}y=\dfrac{w_{c}%
^{2}-w_{b}^{2}+1}{2}. \label{3CartCoordEq}%
\end{equation}

More generally, on a Hilbert space we have:

\begin{theorem}
\label{HilbertCoords}Let $\left(  \mathcal{H},\left\langle \cdot
,\cdot\right\rangle \right)  $ be a real Hilbert space with orthonormal basis
$B.$ The set $C:=B\cup\left\{  \mathbf{0}\right\}  \subset\mathcal{H}$ is a
metric coordinatizing set.
\end{theorem}

\begin{proof}
For $u,v\in\mathcal{H}$ assume $d\left(  u,c\right)  =d\left(  v,c\right)  $
for all $c\in C.$ Then since $\mathbf{0}$ is in $C$ we have $\left\langle
u,u\right\rangle =\left\langle v,v\right\rangle .\ $Further%
\begin{align*}
\left\langle u-c,u-c\right\rangle  &  =\left\langle v-c,v-c\right\rangle \\
\left\langle u,u\right\rangle -2\left\langle c,u\right\rangle +\left\langle
c,c\right\rangle  &  =\left\langle v,v\right\rangle -2\left\langle
c,v\right\rangle +\left\langle c,c\right\rangle \\
\left\langle c,u\right\rangle  &  =\left\langle c,v\right\rangle
\end{align*}
for all $c\in B$ so that $u=v$.
\end{proof}

Using the basis $B$ write an element $w\in\mathcal{H}$ in orthonormal
coordinates as $w\mathbf{=}\left(  \widetilde{w}_{c}\right)  $ where
$\widetilde{w}_{c}=\left\langle w,c\right\rangle $ for each $c\in B.$ Any
point $w\in\mathcal{H}$ is given in metric coordinates by $w=\left(
w_{c}\right)  _{c\in B\cup\left\{  \mathbf{0}\right\}  }$ where $w_{c}%
:=\left\|  w-c\right\|  =d\left(  w,c\right)  .$ With this, the conversion
formulae are%
\begin{align}
\widetilde{w}_{c}  &  =\dfrac{w_{\mathbf{0}}^{2}-w_{c}^{2}+1}{2}\qquad c\in
B\label{HilbertConv}\\
w_{c}  &  =\left(  \left\|  w\right\|  ^{2}-2\widetilde{w}_{c}+1\right)
^{1/2}\qquad c\in B\label{HilbertConv2}\\
w_{\mathbf{0}}  &  =\left\|  w\right\| \nonumber
\end{align}
a straightforward generalization of the finite dimensional formulae, $\left(
\ref{3CartCoordEq}\right)  $ and $\left(  \ref{3metCoordEq}\right)  $.

$\left(  \ref{HilbertConv2}\right)  $ results from the easy calculation%
\begin{align*}
w_{c}  &  =\left\|  w-c\right\|  =\left\langle w-c,w-c\right\rangle ^{1/2}\\
&  =\left(  \left\langle w,w\right\rangle -\left\langle w,c\right\rangle
-\left\langle c,w\right\rangle +\left\langle c,c\right\rangle \right)
^{1/2}\\
&  =\left(  \left\|  w\right\|  ^{2}-2\widetilde{w}_{c}+1\right)  ^{1/2}.
\end{align*}
Solving this equation for $\widetilde{w}_{c}$ yields $\left(
\ref{HilbertConv}\right)  $.

\begin{example}
One must be careful in applying these formulae. They do not necessarily work
on non-Hilbert vector spaces. The finite dimensional Banach space
$\mathbb{R}^{2}$ with the infinity norm has basis $\left\{  \left(
1,0\right)  ,\left(  0,1\right)  \right\}  $ which does not produce a
coordinatizing set in the above manner. Refer to Example \ref{InfinityMetric}.
\end{example}

\section{Derivatives}

A \textbf{curve} in a metric space is a map $\phi:\left(  t_{1},t_{2}\right)
\rightarrow X$ continuous with respect to the metric on $X$ where $\left(
t_{1},t_{2}\right)  $ is a subinterval of $\mathbb{R}$.

\begin{definition}
\label{diffDef}Let $\left(  M,d,X,C\right)  $ be a metric coordinate system
and let $\phi$ be a curve in $X$. Write $\phi$ in metric coordinates as
$\phi_{C}\left(  t\right)  =\left(  \phi_{c}\left(  t\right)  \right)  _{c\in
C}$ . Assuming the limits exist, we define the \textbf{metric-coordinate
derivative }of $\phi$ with respect to $C$ to be $\phi_{C}^{\prime}\left(
t\right)  :=\left(  \phi_{c}^{\prime}\left(  t\right)  \right)  _{c\in C}%
\in\mathbb{R}^{C}$ where%
\[
\phi_{c}^{\prime}\left(  t\right)  =\underset{h\rightarrow0}{\lim}\dfrac
{\phi_{c}\left(  t+h\right)  -\phi_{c}\left(  t\right)  }{h}.
\]
Similarly, the \textbf{forward metric-coordinate derivative} of $\phi$ with
respect to $C$ is $\phi_{C}^{+}\left(  t\right)  :=\left(  \phi_{c}^{+}\left(
t\right)  \right)  _{c\in C}$ where%
\[
\phi_{c}^{+}\left(  t\right)  =\underset{h\rightarrow0^{+}}{\lim}\dfrac
{\phi_{c}\left(  t+h\right)  -\phi_{c}\left(  t\right)  }{h}.
\]

Two curves $\phi$ and $\psi$ are said to be \textbf{(forward)}
\textbf{metric-coordinate-wise} \textbf{tangent} at $t_{0}$ with respect to
$C$ if they meet at $t_{0}$ and have the same (forward) metric-coordinate
derivative, i.e.,%
\begin{align*}
\phi\left(  t_{0}\right)   &  =\psi\left(  t_{0}\right) \\
\phi_{C}^{\prime}\left(  t_{0}\right)   &  =\psi_{C}^{\prime}\left(
t_{0}\right)  \text{\quad(or }\phi_{C}^{+}\left(  t_{0}\right)  =\psi_{C}%
^{+}\left(  t_{0}\right)  \text{).}%
\end{align*}

If there exists $r<\infty$ such that $\left|  \phi_{c}^{\prime}\left(
t_{0}\right)  \right|  \leq r$ for all $c$ then $\phi$ is said to have
\textbf{bounded metric-coordinate speed\footnote{There are other inequivalent
notions of speed such as \textbf{metric speed} $s\left(  t\right)
:=\underset{h\rightarrow0}{\lim}\frac{d\left(  \phi\left(  t+h\right)
,\phi\left(  t\right)  \right)  }{\left|  h\right|  }$ or \textbf{length
speed} $s^{\ast}\left(  t\right)  :=\underset{h\rightarrow0}{\lim
}\frac{L\left(  \phi|_{t}^{t+h}\right)  }{\left|  h\right|  }$ where $L$
refers to the length of the curve.}} at $t_{0}.$
\end{definition}

\begin{remark}
For finite coordinatizing sets $C,$ every metric-coordinate-wise
differentiable curve has bounded metric-coordinate speed at any particular $t
$ in its domain.
\end{remark}

\begin{remark}
\label{NonDiffCoord}A curve $\phi$ in a metric coordinate system which runs
through a coordinatizing point with positive metric speed $s\left(  t\right)
$ can never be differentiable in all of its coordinates; when the curve hits
$c\in C$. The $c$-th coordinate derivative of $\phi$ discontinuously changes
from negative to positive. $\phi$ may still be metric-coordinate-wise
differentiable with respect to another nonintersecting metric coordinatizing
set. Choosing $C$ outside of the region of interest is the reason for the
artifice $\left(  M,d,X,C\right)  $ instead of simply $\left(  X,d,C\right)  $.

Such representational problems are nothing new. E.g., polar coordinates make
do with a continuum-sized discontinuity in representing position.
\end{remark}

The next theorem shows that Definition \ref{diffDef} faithfully generalizes
the traditional derivative on $\mathbb{R}^{n}.$

\begin{theorem}
\label{metCdiffThm}Let $U$ be an open subset of $\mathbb{R}^{n}.$ Let $C$ be a
coordinatizing set for $U$ with respect to the Euclidean metric. Let
$\phi:\left(  t_{1},t_{2}\right)  \rightarrow U$ be a curve and $t\in\left(
t_{1},t_{2}\right)  $ such that $\phi\left(  t\right)  \notin C$. Then $\phi$
is differentiable at $t$ (in the traditional sense) iff it is
metric-coordinate-wise differentiable at $t$.
\end{theorem}

\begin{proof}
First assume $\phi$ is differentiable in the traditional sense at $t.$ Then
for any $c\in C$%
\begin{align*}
\phi_{c}^{\prime}\left(  t\right)   &  =\underset{h\rightarrow0}{\lim}%
\dfrac{\phi_{c}\left(  t+h\right)  -\phi_{c}\left(  t\right)  }{h}\\
&  =\underset{h\rightarrow0}{\lim}\dfrac{d\left(  \phi\left(  t+h\right)
,c\right)  -d\left(  \phi\left(  t\right)  ,c\right)  }{h}\\
&  =\underset{h\rightarrow0}{\lim}\dfrac{f\left(  t+h\right)  -f\left(
t\right)  }{h}%
\end{align*}
where $f\left(  t\right)  =d\left(  \phi\left(  t\right)  ,c\right)  .$ The
Euclidean distance is differentiable except at 0. Since $\phi\left(  t\right)
\neq c$ the function $f$ is the composition of two differentiable functions
and hence differentiable. Thus $\phi$ is metric-coordinate-wise differentiable
at $t.$

The converse is slightly more difficult. Assume $\phi$ is
metric-coordinate-wise differentiable at $t$ with respect to $C$. We prove
that $\phi$ is differentiable in the traditional sense in the context of
$\mathbb{R}^{2};$ the generalization to $\mathbb{R}^{n}$ is immediate. There
exist two points from $C$, say $a$ and $b\in\mathbb{R}^{2}$, which together
with $\phi\left(  t\right)  \in\mathbb{R}^{2}$ are non-colinear--else $C$
would not effectively discriminate between all points of $U$. Define
$\binom{f}{g}:\mathbb{R}^{2}\rightarrow\mathbb{R}^{2}$ by $f\left(  x\right)
:=\left\|  x-a\right\|  $, and $g\left(  x\right)  :=\left\|  x-b\right\|  $
so that $\binom{f}{g}$ is differentiable at any point not equal to $a$ or $b.$
We will show that $\binom{f}{g}^{-1}$ and $\binom{f}{g}\circ\phi$ are
differentiable so that $\phi$ is the composition of two differentiable
functions. Our assumption that $\phi$ is coordinate-wise differentiable at $t$
immediately gives the differentiability of $\binom{f}{g}\circ\phi:\left(
t_{1},t_{2}\right)  \rightarrow\mathbb{U}$ at $t$. To prove the
differentiability of $\binom{f}{g}^{-1},$ we show $D\binom{f}{g}$ is
nonsingular at $x=\phi\left(  t\right)  $:%
\begin{align*}
D\binom{f}{g}\left(  x\right)   &  =\left[
\begin{array}
[c]{cc}%
f_{u}\left(  x\right)  & g_{u}\left(  x\right) \\
f_{v}\left(  x\right)  & g_{v}\left(  x\right)
\end{array}
\right] \\
&  =\left[
\begin{array}
[c]{cc}%
\dfrac{x_{1}-a_{1}}{\left\|  x-a\right\|  } & \dfrac{x_{1}-b_{1}}{\left\|
x-b\right\|  }\\
\dfrac{x_{2}-a_{2}}{\left\|  x-a\right\|  } & \dfrac{x_{2}-b_{2}}{\left\|
x-b\right\|  }%
\end{array}
\right]
\end{align*}
which is singular if and only if one column vector is a multiple of the other,
i.e.,
\[
\left(  \dfrac{x_{1}-a_{1}}{\left\|  x-a\right\|  },\dfrac{x_{2}-a_{2}%
}{\left\|  x-a\right\|  }\right)  =\lambda\left(  \dfrac{x_{1}-b_{1}}{\left\|
x-b\right\|  },\dfrac{x_{2}-b_{2}}{\left\|  x-b\right\|  }\right)
\]
or equivalently that
\begin{align*}
0  &  =\left(  x_{1}-a_{1},x_{2}-a_{2}\right)  -\lambda_{1}\left(  x_{1}%
-b_{1},x_{2}-b_{2}\right) \\
&  =\left(  1-\lambda_{1}\right)  x-a-\lambda_{1}b.
\end{align*}
In this case $x,$ $a$ and $b$ are colinear, which cannot happen. Thus
$D\binom{f}{g}$ is nonsingular, $\binom{f}{g}$ is locally invertible, and
$\binom{f}{g}^{-1}$ is differentiable. Hence $\binom{f}{g}^{-1}\circ\binom
{f}{g}\circ\phi=\phi$ is differentiable.
\end{proof}

In view of this theorem we could use ``differentiable'' in lieu of the awkward
phrasing ``metric-coordinate-wise differentiable''. But in order to be
perfectly clear in this nascent setting we usually employ the full term.

\begin{example}
Any curve in $\left(  \mathbb{E}^{2},d,H^{2},\left\{  a,b\right\}  \right)  $
from Example \ref{H^2CoordEx} which satisfies the conditions of Proposition
\ref{metTriangleProp}%
\begin{align*}
\phi_{a},\phi_{b}  &  \geq0\\
\left|  \phi_{a}\left(  t\right)  -\phi_{b}\left(  t\right)  \right|   &
\leq1\\
\phi_{a}\left(  t\right)  +\phi_{b}\left(  t\right)   &  \geq1,
\end{align*}
and is differentiable in each of its coordinates will be a
metric-coordinate-wise differentiable curve$.$
\end{example}

\begin{example}
In the Hilbert space $L^{2}\left(  \mathbb{R}\right)  $ define the curve
$\phi:\left(  -\infty,\infty\right)  \rightarrow L^{2},$ by $\phi\left(
t\right)  \left(  x\right)  :=\chi_{\left[  0,1\right]  }\left(  x-t\right)  $
where $\chi_{\left[  0,1\right]  }$ is the characteristic function of the unit
interval. $\phi$ is not Frechet differentiable, but it is
metric-coordinate-wise differentiable with respect to certain metric
coordinatizing sets as we show.

The difference quotient
\[
\frac{\phi\left(  t+h\right)  -\phi\left(  t\right)  }{h}%
\]
does not converge in $L^{2}$ to any $g\in L^{2};$ it does however converge in
the distribution sense to the difference of Dirac deltas $\delta\left(
t+1\right)  -\delta\left(  t\right)  .$

Choose an orthonormal basis $B$ of $L^{2}$ consisting of continuous functions.
A metric coordinate system is then automatically given by $C:=B\cup\left\{
\mathbf{0}\right\}  $ by Theorem \ref{HilbertCoords}. Then%
\begin{gather}
\phi_{c}^{\prime}\left(  t\right)  =\underset{h\rightarrow0}{\lim}%
\frac{\phi_{c}\left(  t+h\right)  -\phi_{c}\left(  t\right)  }{h}%
=\underset{h\rightarrow0}{\lim}\frac{\left\|  \phi\left(  t+h\right)
-c\right\|  _{2}-\left\|  \phi\left(  t\right)  -c\right\|  _{2}}%
{h}\nonumber\\
=\underset{h\rightarrow0}{\lim}\frac{\left\|  \phi\left(  t+h\right)
-c\right\|  _{2}^{2}-\left\|  \phi\left(  t\right)  -c\right\|  _{2}^{2}%
}{h\left(  \left\|  \phi\left(  t+h\right)  -c\right\|  _{2}+\left\|
\phi\left(  t\right)  -c\right\|  _{2}\right)  }\\
=\underset{h\rightarrow0}{\lim}\frac{\int\left(  \left[  \phi\left(
t+h\right)  ^{2}-\phi\left(  t\right)  ^{2}\right]  -2c\left[  \phi\left(
t+h\right)  -\phi\left(  t\right)  \right]  \right)  }{h\left(  \left\|
\phi\left(  t+h\right)  -c\right\|  _{2}+\left\|  \phi\left(  t\right)
-c\right\|  _{2}\right)  } \label{L^2Ex}%
\end{gather}
and since $\phi\left(  t\right)  ^{2}=\phi\left(  t\right)  $ by the nature of
the characteristic function, we have
\begin{gather}
=\underset{h\rightarrow0}{\lim}\frac{\int\left(  \left[  \phi\left(
t+h\right)  -\phi\left(  t\right)  \right]  -2c\left[  \phi\left(  t+h\right)
-\phi\left(  t\right)  \right]  \right)  }{h\left(  \left\|  \phi\left(
t+h\right)  -c\right\|  _{2}+\left\|  \phi\left(  t\right)  -c\right\|
_{2}\right)  }\\
=\underset{h\rightarrow0}{\lim}\frac{\int\left(  \left[  1-2c\right]  \left[
\phi\left(  t+h\right)  -\phi\left(  t\right)  \right]  \right)  /h}{\left(
\left\|  \phi\left(  t+h\right)  -c\right\|  _{2}+\left\|  \phi\left(
t\right)  -c\right\|  _{2}\right)  }\\
=-\frac{2\left[  c\left(  t+1\right)  -c\left(  t\right)  \right]  }{2\left\|
\phi\left(  t\right)  -c\right\|  _{2}}=\frac{c\left(  t\right)  -c\left(
t+1\right)  }{\phi_{c}\left(  t\right)  }. \label{L^2Ex3}%
\end{gather}
The second to last equality in line $\left(  \ref{L^2Ex3}\right)  $ is the
reason we require the continuous basis.
\end{example}

\begin{example}
[Observer dependence of smoothness]\label{CoordDep}On general metric spaces,
metric-coordinate-wise differentiability is crucially dependent on the
particular coordinate system. E.g., let $M:=\mathbb{R}^{2}$ with Euclidean
metric $d.\ $%
\[
X:=\left\{  \left(  x,\left|  x\right|  \right)  :\left|  x\right|
\leq1\right\}  \text{.}%
\]
Two different metric coordinate systems for $X$ are given by the singletons
$C_{1}=\left\{  \left(  -2,0\right)  \right\}  $ and $C_{2}=\left\{  \left(
1,1\right)  \right\}  $. The curve $\psi:\left(  -1,1\right)  \rightarrow X$
given by $\psi\left(  t\right)  :=\left(  t,\left|  t\right|  \right)  $ is
metric-coordinate-wise differentiable at $t=0$ with respect to $\left\{
\left(  -2,0\right)  \right\}  $, but not with respect to $\left\{  \left(
1,1\right)  \right\}  $. I.e., an observer at $\left(  1,1\right)  $ measures
the jarring difference in distance at time $t=0,$ whereas an observer at
$\left(  -2,0\right)  $ measures a smoothly changing distance.

One could give a more involved definition of metric-coordinate-wise
differentiability that eliminates coordinate dependence, but we will not
pursue it here.
\end{example}

\section{Vector fields}

A map $f:\left(  X,d_{X}\right)  \rightarrow\left(  Y,d_{Y}\right)  $ from one
metric space to another is called $K$-\textbf{Lipschitz} (or just Lipschitz)
if
\[
d_{Y}\left(  f\left(  x\right)  ,f\left(  y\right)  \right)  \leq
Kd_{X}\left(  x,y\right)
\]
for all $x,y\in X.$ The map $f$ is called \textbf{locally Lipschitz} if for
each point there is a $K\geq0$ and a neighborhood on which $f$ is $K$-Lipschitz.

Consider the autonomous ordinary differential equation%
\begin{equation}
\dot{x}=f\left(  x\right)  \label{ODE}%
\end{equation}
on a Banach space $B$. $f$ is called the vector field associated with the
differential equation $\left(  \ref{ODE}\right)  $ and is a map
$f:B\rightarrow B.$ The Cauchy-Lipschitz Theorem on Banach spaces guarantees
that if $f$ is locally Lipschitz then unique solutions exist for short time
from any initial condition $x_{0}\in B$. I.e., there exists $x:\left(
-\delta,\delta\right)  \rightarrow B$ for some $\delta>0$ with $x\left(
0\right)  =x_{0}$ satisfying $\left(  \ref{ODE}\right)  .$ The goal of this
section is to achieve a similar result for metric coordinate systems using the
fact that $X$ may be associated with a subset of $\mathbb{R}^{C}$ via the $C$
embedding\footnote{Our focus in this paper is on autonomous dynamics--i.e.,
vector fields which do not change in time--but the mechanics of generating
time-dependent flows may be extended with little extra effort with the
standard trick. Simply work on the metric space $X\times\mathbb{R}$ with
$\mathbb{R}$ representing the time coordinate, then carefully project
solutions on $X.$}.

In order to achieve this goal we use a new metric $d_{C}$ on $X$. We will see
that in many important cases $\left(  X,d_{C}\right)  $ is homeomorphic to
$\left(  X,d\right)  $. For a metric coordinate system $\left(
M,d,X,C\right)  $ define $d_{C}:X\times X\rightarrow\mathbb{R}$ by%
\begin{equation}
d_{C}\left(  x,y\right)  :=\underset{c\in C}{\sup}\left|  x_{c}-y_{c}\right|
\label{dC}%
\end{equation}
for $x,y\in X.$ To see that this gives a finite number for arbitrary
coordinatizing sets $C$ notice that
\begin{equation}
\left|  x_{c}-y_{c}\right|  =\left|  d\left(  x,c\right)  -d\left(
y,c\right)  \right|  \leq d\left(  x,y\right)  \label{dc<d}%
\end{equation}
by the triangle inequality. This shows that $d_{C}\leq d$. Further, a subset
of $X$ is bounded with respect to $d$ if and only if it is bounded with
respect to $d_{C}$.

\begin{definition}
Let $\left(  M,d,X,C\right)  $ be a metric coordinate system. Let $X^{+}$
represent the space of curves $\phi:\left[  0,\delta\right)  \rightarrow
\left(  X,d_{C}\right)  $ which are forward metric-coordinate-wise
differentiable with bounded metric-coordinate-wise speed at $t=0$ and define
an equivalence relation $\sim$ on $X^{+}$ by $\phi\sim\psi$ if $\phi$ is
forward tangent to $\psi$ at $t=0$. The space of equivalence classes for
$\sim$ is the \textbf{tangent bundle }of $X$ and is written with the symbol $TX.$

The set of equivalence classes of curves under $\sim$ for which $\phi\left(
0\right)  =x\in X$ is the \textbf{tangent space} of $X$ at $x$ and is referred
to with the symbol $T_{x}X.$

We also define a metric $d_{C}^{T}$ on the tangent bundle $TX$ by%
\[
d_{C}^{T}\left(  \left[  \phi\right]  ,\left[  \psi\right]  \right)
:=\max\left\{  d_{C}\left(  \phi\left(  0\right)  ,\psi\left(  0\right)
\right)  ,\underset{c\in C}{\sup}\left|  \phi_{c}^{+}\left(  0\right)
-\psi_{c}^{+}\left(  0\right)  \right|  \right\}
\]
where $\phi\in\left[  \phi\right]  \in TX$ and $\psi\in\left[  \psi\right]
\in TX$.
\end{definition}

Forward derivatives are used because of the abundance of metric spaces with
boundaries. Henceforth we only consider forward derivative, but everything
could be formulated in terms of two-sided derivatives as well.

Clearly $TX$ is the disjoint union $\underset{x\in X}{\amalg}T_{x}X$. Notice
that $TX$ depends on $C$, not just $X$ and that we use the metric $d_{C}$
instead of $d$.

\begin{remark}
Any member $\left[  v\right]  \in T_{x}X$ is by definition an equivalence
class of curves, but may be represented with a single element of
$\mathbb{R}^{C}.$ This is true since any two members $v,w\in\left[  v\right]
$ have the same forward metric-coordinate derivatives at 0, i.e.,%
\[
v_{C}^{+}\left(  0\right)  =\left(  v_{c}^{+}\left(  0\right)  \right)  _{c\in
C}=\left(  w_{c}^{+}\left(  0\right)  \right)  _{c\in C}=w_{C}^{+}\left(
0\right)  \in\mathbb{R}^{C}.
\]
It would not be too egregious an abuse of notation to write $T_{x}%
X\subset\mathbb{R}^{C}.$
\end{remark}

\begin{remark}
Though the symbol $T_{x}X$ represents a vector space in the context of
differentiable manifolds, this is often not true in metric coordinate systems.
E.g., from Example \ref{H^2CoordEx} we consider the closed half space
$\overline{H^{2}}$ and the metric coordinate system $\left(  \mathbb{E}%
^{2},d,\overline{H^{2}},\left\{  a,b\right\}  \right)  $. Then $T_{x}%
\overline{H^{2}}\ $is naturally identified with $\mathbb{R}^{2},$ a vector
space, for any interior point $x$ of $\overline{H^{2}}$, but this is not true
for $x$ on the boundary. Notice however that for any system $\left(
M,d,X,C\right)  $ the tangent space $T_{x}X$ consists of rays emanating from
the origin, since curves may be reparametrized to have greater or smaller
metric coordinate derivative.
\end{remark}

\begin{definition}
On a metric coordinate system $\left(  M,d,X,C\right)  $ a
\textbf{metric-coordinate vector field} is a map $V:X\rightarrow TX$ such that
$V\left(  x\right)  \in T_{x}X$ for each $x\in X$ with $V\left(  x\right)
_{c}$ uniformly bounded in $c$ for each $x,$ i.e., $V\left(  x\right)  $ has
bounded metric coordinate speed for each $x$.

Such a vector field is called (\textbf{locally}) \textbf{Lipschitz} if
$V:\left(  X,d_{C}\right)  \rightarrow\left(  TX,d_{C}^{T}\right)  $ is
(locally) Lipschitz.

A \textbf{solution} to $V$ with \textbf{initial condition} $x\in X$ is a curve
$\sigma:\left[  0,\delta\right)  \rightarrow\left(  X,d_{C}\right)  $ for some
$\delta>0$ with $\sigma\left(  0\right)  =x$ such that $\sigma^{+}\left(
t\right)  =V\left(  \sigma\left(  t\right)  \right)  $ for all $t\in\left[
0,\delta\right)  $.

A metric-coordinate vector field on $\left(  M,d,X,C\right)  $ is said to have
\textbf{unique solutions }if for any point $x\in X,$ there exists a solution
$\sigma:\left[  0,\delta\right)  \rightarrow X$ with $\sigma\left(  0\right)
=x,$ and if $\tau:\left[  0,\epsilon\right)  \rightarrow X$ is another
solution with $\tau\left(  0\right)  =x,$ then for $t\in\left[  0,\min\left\{
\delta,\epsilon\right\}  \right]  $ we have $\tau\left(  t\right)
=\sigma\left(  t\right)  $.
\end{definition}

\begin{remark}
Any solution to a metric-coordinate vector field with unique solutions may be
continued to produce a solution with maximal domain using a straightforward
analytic argument.
\end{remark}

\section{\label{CauLipSection}Cauchy-Lipschitz Theorem for metric coordinate systems}

\begin{theorem}
\label{CauLips}Let $\left(  M,d,X,C\right)  $ be a metric coordinate system
and assume $\left(  X,d_{C}\right)  $ is locally compact. Let $V:\left(
X,d_{C}\right)  \rightarrow\left(  TX,d_{C}^{T}\right)  $ be a locally
Lipschitz metric-coordinate vector field. Then $V$ has unique solutions.
\end{theorem}

This section is devoted to the proof. Much of the following could be
conceptually simplified by considering only finite metric coordinatizing sets.
But the setting of a metric space is so abstract that it is a great advantage
to consider arbitrary $C$.

The outline of the proof begins by viewing $X$ and $TX$ as subsets of
$\mathbb{R}^{C}.$ We then extend the vector field $V$ to a map $V^{1}%
:\mathbb{R}^{C}\rightarrow\mathbb{R}^{C},$ use the traditional
Cauchy-Lipschitz Theorem to guarantee solutions, and verify that restrictions
of these solutions to $X$ remain in $X$ for short time with the Nagumo
Invariance Theorem.

The problem with this plan is that $\mathbb{R}^{C}$ with the supremum norm is
not a Banach space when $C$ is infinite, and so the standard Cauchy-Lipschitz
theorem (Theorem \ref{VSolRn} below) does not apply. However, the space of
bounded $C$-tuples%
\[
\mathbb{R}_{b}^{C}:=\left\{  x_{C}\in\mathbb{R}^{C}:\underset{c\in C}{\sup
}\left|  x_{c}\right|  <\infty\right\}
\]
is a Banach space for any set $C$ with norm\footnote{We reserve the notation
$\left\|  \cdot\right\|  $ for this supremum norm henceforth.}%
\[
\left\|  x\right\|  :=\underset{c\in C}{\sup}\left|  x_{c}\right|  \text{.}%
\]
To carry out our plan, we embed $X$ into $\mathbb{R}_{b}^{C}$ instead.

Let $w\in X$ be a distinguished element (arbitrarily chosen), and for each
$x\in X$ define the embedding $i:\left(  X,d_{C}\right)  \rightarrow\left(
\mathbb{R}_{b}^{C},\left\|  \cdot\right\|  \right)  $ by%
\[
i\left(  x\right)  _{c}:=x_{c}-d\left(  c,w\right)  .
\]
Then $i\left(  x\right)  \in\mathbb{R}_{b}^{C}$ since%
\[
i\left(  x\right)  _{c}=d\left(  x,c\right)  -d\left(  c,w\right)  \leq
d\left(  x,w\right)
\]
which is uniformly bounded in $C$. Subtracting $d\left(  c,w\right)  $ in the
definition of $i$ is only necessary in the case that $C$ is unbounded in the
metric sense. Finally $i$ is an isometry (in particular it is injective) since%
\begin{gather*}
\left\|  i\left(  x\right)  -i\left(  y\right)  \right\|  =\underset{c\in
C}{\sup}\left|  i\left(  x\right)  _{c}-i\left(  y\right)  _{c}\right| \\
=\underset{c\in C}{\sup}\left|  d\left(  x,c\right)  -d\left(  c,w\right)
-\left[  d\left(  y,c\right)  -d\left(  c,w\right)  \right]  \right| \\
=\underset{c\in C}{\sup}\left|  d\left(  x,c\right)  -d\left(  y,c\right)
\right|  =d_{C}\left(  x,y\right)  .
\end{gather*}

We will need the following results.

\begin{theorem}
[Cauchy-Lipschitz]\label{VSolRn}A locally Lipschitz vector field on a Banach
space has unique solutions.
\end{theorem}

Here we are referring to the traditional notion of vector field, not
metric-coordinate vector fields. Proofs are legion.

\begin{remark}
The uniqueness of one-sided solutions, required for this section, is also
true. See \cite{CalcNBleecker}, e.g.
\end{remark}

\begin{theorem}
[Lipschitz Extension]\label{VFExtension}If $S$ is a subset of a metric space
$(X,d),$ and if $f:S\rightarrow\mathbb{R}$ is $K$-Lipschitz, then
$\overline{f}:X\rightarrow\mathbb{R}$ defined by%
\[
\overline{f}(x):=\sup\left\{  f\left(  y\right)  -K\cdot d\left(  x,y\right)
|y\in S\right\}
\]
equals $f$ on $S$ and is $K$-Lipschitz.
\end{theorem}

\begin{proof}
Given in $\cite{McShane}$.
\end{proof}

\begin{lemma}
\label{LipExtLemma}If $S$ is a subset of a metric space $\left(  X,d\right)
,$ and if $f:S\rightarrow\mathbb{R}_{b}^{C}$ is $K$-Lipschitz, then there
exists a $K$-Lipschitz extension $\overline{f}:\left(  X,d\right)
\rightarrow\mathbb{R}_{b}^{C}.$
\end{lemma}

\begin{proof}
Use the Lipschitz Extension Theorem on each coordinate to get $\overline
{f}:X\rightarrow\mathbb{R}^{C}$ which is $K$-Lipschitz in each coordinate. We
need to check that $\overline{f}\left(  X\right)  \subset\mathbb{R}_{b}^{C}.$
For any $x\in X$ and $y\in S$
\begin{align*}
\underset{c\in C}{\sup}\left|  \overline{f}_{c}\left(  x\right)  \right|   &
\leq\underset{c\in C}{\sup}\left\{  \left|  \overline{f}_{c}\left(  x\right)
-\overline{f}_{c}\left(  y\right)  \right|  +\left|  \overline{f}_{c}\left(
y\right)  \right|  \right\} \\
&  \leq Kd\left(  x,y\right)  +\left\|  f\left(  y\right)  \right\|  <\infty.
\end{align*}
Therefore $\overline{f}_{C}\left(  x\right)  \in\mathbb{R}_{b}^{C}$.
\end{proof}

The \textbf{upper forward derivative} of a function $f:\left[  a,b\right]
\rightarrow\mathbb{R}$ is defined by%
\[
D^{+}f\left(  t\right)  :=\underset{h\rightarrow0^{+}}{\overline{\lim}%
}\frac{f\left(  t+h\right)  -f\left(  t\right)  }{h}\text{.}%
\]

\begin{lemma}
Let $f:\left[  a,b\right]  \rightarrow\mathbb{R}$ be continuous with
$D^{+}f\left(  t\right)  \leq Kf\left(  t\right)  .$ Then $f\left(  t\right)
\leq f\left(  a\right)  $ for all $t\in\left[  a,b\right]  .$
\end{lemma}

\begin{proof}
See \cite[p. 354]{Titchmarsh} for the following result: Let $F:\left[
a,b\right]  \rightarrow\mathbb{R}$ be continuous with $D^{+}F\leq0$ for all
$t\in\left[  a,b\right)  .$ Then $F\left(  a\right)  \geq F\left(  b\right)  .$

Now apply this to $F\left(  t\right)  :=e^{-Kt}f\left(  t\right)  $.
\end{proof}

\begin{definition}
A subset $S$ of a normed vector space $E$ is said to be \textbf{positively
invariant} with respect to the vector field $V:E\rightarrow E$ if any forward
solution $\sigma:\left[  0,\delta\right)  \rightarrow E$ to $V$ with initial
condition $\sigma\left(  0\right)  \in S$ has $\sigma\left(  t\right)  \in S$
for all $t\in\left[  0,\delta\right)  .$
\end{definition}

For a point $x$ in a metric space $\left(  X,d\right)  $ and a subset $S,$ the
distance from $x$ to $S$ is defined as%
\[
d\left(  x,S\right)  :=\underset{y\in S}{\inf}\left\{  d\left(  x,y\right)
\right\}  =:d\left(  S,x\right)  .
\]
It is easy to check that $d\left(  x,S\right)  \leq d\left(  x,y\right)
+d\left(  y,S\right)  $ for any $y\in X.$ As a consequence the distance is
continuous in $x$.

\begin{theorem}
[Nagumo Invariance]Let $E$ be a normed vector space space, let $V:E\rightarrow
E$ be a map, and let $S$ be a closed subset of $E$. Suppose that at each $a\in
S$ the vector field $V$ is tangent to $S$ in the following sense: there exists
an open neighborhood $\Omega_{a}$ and $K_{a}>0$ such that%
\begin{equation}
\underset{h\rightarrow0^{+}}{\overline{\lim}}\frac{d\left(  x+hV\left(
x\right)  ,S\right)  -d\left(  x,S\right)  }{h}\leq K_{a}d\left(  x,S\right)
\label{VolkCondition}%
\end{equation}
for all $x\in\Omega_{a},$ where $d$ is the metric induced by the norm on $E.$

Then $S$ is positively invariant with respect to the vector field $V.$
\end{theorem}

\begin{proof}
This generalization of Nagumo's result on $\mathbb{R}^{n}$ is due to Volkmann
and is given in \cite{Volkmann} under more general conditions. We adapt his
proof to this context. Similar results are surveyed in \cite[pp.
70-71,98]{Pavel}.

Assume $S$ is not positively invariant. Then there is a solution
$\sigma:\left[  0,\delta\right)  \rightarrow E$ with $\sigma\left(  0\right)
\in S$ and $\sigma\left(  t_{0}\right)  \notin S$ for some $t_{0}\in\left[
0,\delta\right)  .$ Let
\[
t_{1}:=\sup\left\{  t:\sigma\left(  \left[  0,t\right)  \right)  \subset
S\right\}  .
\]
Since $S$ is closed, $a=\sigma\left(  t_{1}\right)  \in S$ and $0\leq
t_{1}<t_{0}<\delta$.

For the point $\sigma\left(  t_{1}\right)  $ choose $\Omega$ according to the
assumptions of the theorem and let $t_{2}$ be chosen greater than $t_{1}$ such
that $\sigma\left(  t\right)  \in\Omega$ for $t\in\left[  t_{1},t_{2}\right]
$. By the definition of $t_{1}$ there exists some $t_{3}\in\left(  t_{1}%
,t_{2}\right)  $ such that $\sigma\left(  t_{3}\right)  \notin S.$

Define $\eta:\left[  t_{1},t_{2}\right)  \rightarrow\left[  0,\infty\right)  $
by%
\[
\eta\left(  s\right)  :=d\left(  \smallskip\sigma\left(  s\right)  ,S\right)
.
\]
Certainly $\eta$ is continuous, positive, and $\eta\left(  t_{1}\right)  =0.$
We prove that the upper forward derivative of $\eta$ is less than $K_{a}\eta$
on its domain, so that $\eta\left(  s\right)  \equiv0$ by the previous lemma.
To this end, fix $s\in\left[  t_{1},t_{2}\right)  .$ Then for $h>0$%
\begin{align*}
\eta\left(  s+h\right)   &  =d\left(  \smallskip\sigma\left(  s+h\right)
,S\right) \\
&  \leq d\left(  \smallskip\sigma\left(  s+h\right)  ,\sigma\left(  s\right)
+hV\left(  \sigma\left(  s\right)  \right)  \right)  +d\left(  \smallskip
\sigma\left(  s\right)  +hV\left(  \sigma\left(  s\right)  \right)  ,S\right)
\\
&  =\left\|  \sigma\left(  s+h\right)  -\sigma\left(  s\right)  -hV\left(
\sigma\left(  s\right)  \right)  \right\|  +d\left(  \smallskip\sigma\left(
s\right)  +hV\left(  \sigma\left(  s\right)  \right)  ,S\right) \\
&  =\text{o}\left(  h\right)  +d\left(  \smallskip\sigma\left(  s\right)
+hV\left(  \sigma\left(  s\right)  \right)  ,S\right)
\end{align*}

The last equality results from the fact that $\sigma$ is a solution to $V.$
Thus the upper forward derivative of $\eta$ is%
\begin{align*}
\overline{D}_{\eta}^{+}\left(  s\right)   &  :=\underset{h\rightarrow0^{+}%
}{\overline{\lim}}\frac{\eta\left(  s+h\right)  -\eta\left(  s\right)  }{h}\\
&  \leq\underset{h\rightarrow0^{+}}{\overline{\lim}}\frac{d\left(
\smallskip\sigma\left(  s\right)  +hV\left(  \sigma\left(  s\right)  \right)
,S\right)  -d\left(  \smallskip\sigma\left(  s\right)  ,S\right)  }{h}\\
&  \leq K_{a}d\left(  \smallskip\sigma\left(  s\right)  ,S\right)  =K_{a}%
\eta\left(  s\right)  .
\end{align*}
The last inequality is from $\left(  \ref{VolkCondition}\right)  .$ Thus
$\eta\left(  s\right)  \equiv0$ so that $\sigma\left(  t\right)  \in S$ for
all $t\in\left[  t_{1},t_{2}\right)  ,$ contradicting $\sigma\left(
t_{3}\right)  \notin S.$
\end{proof}

Finally we are ready to prove the major result.

\begin{proof}
[Proof of Theorem \ref{CauLips}]Let $x_{0}\in X.$ Since $\left(
X,d_{C}\right)  $ is locally compact, there exists a compact ball
$\overline{B}:=\overline{B}_{d_{C}}\left(  x_{0},r\right)  $ for some $r>0$.
We may assume $r$ is chosen small enough so that $V$ is $K$-Lipschitz on
$\overline{B}.$ Notice that the imbedding map $i $ gives $i\left(
\overline{B}\right)  \subset\overline{B}_{\left\|  \cdot\right\|  }\left(
i\left(  x_{0}\right)  ,r\right)  $ where $B_{d_{C}}$ refers to a ball in
$\left(  X,d_{C}\right)  $ and $B_{\left\|  \cdot\right\|  }$ refers to a ball
in $\mathbb{R}_{b}^{C}.$ Further $i\left(  \overline{B}\right)  $ is compact,
being the continuous (isometric) image of the compact space $\overline{B}.$

The metric-coordinate vector field $V:X\rightarrow TX$ transfers to a map
$V^{1}$ on $i\left(  B\right)  $ via the following diagram:%
\[%
\begin{array}
[c]{ccc}%
\qquad V:\left(  \overline{B},d_{C}\right)  & \overset{K\text{-Lip}%
}{\rightarrow} & \left(  TX,d_{C}^{T}\right)  \qquad\qquad\qquad\\
\text{{\small (isometry) }}i\text{ }\downarrow\qquad &  & \qquad
\downarrow\text{ }\pi\text{{\small (weak contraction)}}\\
\qquad V^{1}:\left(  i\left(  \overline{B}\right)  ,\left\|  \cdot\right\|
\right)  & \rightarrow & \left(  \mathbb{R}_{b}^{C},\left\|  \cdot\right\|
\right)  \qquad\qquad\qquad
\end{array}
\]
where $\pi\left(  \left[  \phi\right]  \right)  :=\phi_{C}^{+}\left(
0\right)  .$ The map $\pi$ is a weak contraction (i.e., $K$-Lipschitz with
$K\leq1$) since%
\[
\left\|  \pi\left(  \left[  \phi\right]  \right)  -\pi\left(  \left[
\psi\right]  \right)  \right\|  =\left\|  \phi_{C}^{+}\left(  0\right)
-\psi_{C}^{+}\left(  0\right)  \right\|  \leq d_{C}^{T}\left(  \left[
\phi\right]  ,\left[  \psi\right]  \right)  .
\]
Notice $V^{1}\circ i=\pi\circ V$ so we see that $V^{1}$ is $K$-Lipschitz
since
\begin{gather*}
\left\|  V^{1}\left(  i\left(  x\right)  \right)  -V^{1}\left(  i\left(
y\right)  \right)  \right\|  =\left\|  \pi\left(  V\left(  x\right)  \right)
-\pi\left(  V\left(  y\right)  \right)  \right\| \\
\leq d_{C}^{T}\left(  V\left(  x\right)  ,V\left(  y\right)  \right)  \leq
Kd_{C}\left(  x,y\right)  =K\left\|  i\left(  x\right)  -i\left(  y\right)
\right\|  .
\end{gather*}

Extend $V^{1}$ to a Lipschitz vector field $V^{2}$ on all of $\mathbb{R}%
_{b}^{C}$ via Lemma \ref{LipExtLemma}. We will prove that a solution to
$V^{2}$ starting at $i\left(  x_{0}\right)  $ remains in $i\left(
\overline{B}\right)  $ for short time. Modify $V^{2}$ to be an invariant
vector field on $B_{\left\|  \cdot\right\|  }\left(  i\left(  x_{0}\right)
,r\right)  $ by shrinking the speed to 0 near its boundary. To do this define
the new vector field $V^{3}:\mathbb{R}_{b}^{C}\rightarrow\mathbb{R}_{b}^{C}$
to be%
\[
V^{3}\left(  w\right)  :=\left\{
\begin{array}
[c]{ll}%
V^{2}\left(  w\right)  & w\in B_{\left\|  \cdot\right\|  }\left(  i\left(
x_{0}\right)  ,r/2\right) \\
0 & w\notin B_{\left\|  \cdot\right\|  }\left(  i\left(  x_{0}\right)
,r\right) \\
\left(  2-\frac{2}{r}\left\|  w-i\left(  x_{0}\right)  \right\|  \right)
V^{2}\left(  w\right)  \text{ } & \text{ otherwise}%
\end{array}
\right.
\]
which is again Lipschitz (which is verified in Lemma \ref{fLipExten} below),
say with constant $K_{1}.$ The Cauchy-Lipschitz Theorem on Banach spaces then
provides unique solutions to $V^{3}.$

For the penultimate step of the proof we invoke the Nagumo Invariance Theorem
to demonstrate that the solutions to $V^{3}$ which begin in $i\left(
\overline{B}\right)  $ remain in $i\left(  \overline{B}\right)  .$ We use
$\Omega=\mathbb{R}_{b}^{C}.$ We use the metric $d_{\infty}$ derived from the
norm $\left\|  \cdot\right\|  $ on $\mathbb{R}_{b}^{C}$. First consider $w\in
i\left(  \overline{B}\right)  ;$ we get%
\begin{align}
&  \frac{d_{\infty}\left(  \smallskip w+hV^{3}\left(  w\right)  ,i\left(
\overline{B}\right)  \right)  -d_{\infty}\left(  \smallskip w,i\left(
\overline{B}\right)  \right)  }{h}\nonumber\\
&  =\frac{d_{\infty}\left(  \smallskip w+hV^{3}\left(  w\right)  ,i\left(
\overline{B}\right)  \right)  }{h}\nonumber\\
&  \leq\frac{d_{\infty}\left(  \smallskip w+hV^{3}\left(  w\right)
,\phi\left(  h\right)  \right)  +d_{\infty}\left(  \smallskip\phi\left(
h\right)  ,i\left(  \overline{B}\right)  \right)  }{h} \label{c(h)Estim1}%
\end{align}
where $\phi:\left[  0,\delta\right)  \rightarrow i\left(  \overline{B}\right)
$ is a curve with $\phi\left(  0\right)  =w$ and $\phi^{+}\left(  0\right)
=V^{3}\left(  w\right)  .$ It is not immediately clear that there is such a
curve which remains in $i\left(  \overline{B}\right)  $. To see that such a
$\phi$ exists consider the three cases:

1. If $w\in B_{\left\|  \cdot\right\|  }\left(  i\left(  x_{0}\right)
,r/2\right)  $ then $V^{3}\left(  w\right)  =\pi\left(  V\left(  i^{-1}\left(
w\right)  \right)  \right)  $ so that there exists a member of the equivalence
class $V\left(  i^{-1}\left(  w\right)  \right)  ,$ call it $\psi:\left[
0,\delta\right)  \rightarrow X$ with $\psi^{+}\left(  0\right)  =\pi\left(
V\left(  i^{-1}\left(  w\right)  \right)  \right)  .$ We assume $\delta>0$ is
chosen small enough that $\psi$ remains in $B_{d_{C}}\left(  x_{0},r/2\right)
$ which may be done since $\psi$ is continuous with respect to $d_{C}$. Then
$\phi:=i\circ\psi$ is the desired curve.

2. If $w\in B_{\left\|  \cdot\right\|  }\left(  i\left(  x_{0}\right)
,r\right)  \backslash B_{\left\|  \cdot\right\|  }\left(  i\left(
x_{0}\right)  ,r/2\right)  ,$ the same approach as Case 1 works again; just
reparametrize with multiplicative factor
\[
\left(  2-\frac{2}{r}\left\|  w-i\left(  x_{0}\right)  \right\|  \right)  .
\]

3. If $w\notin B_{\left\|  \cdot\right\|  }\left(  i\left(  x_{0}\right)
,r\right)  $ use the constant curve $\phi\left(  t\right)  \equiv w.$ This
seems simple, but it is the reason we modified $V^{2}$ to $V^{3}$; when $w$ is
on the boundary of $B_{\left\|  \cdot\right\|  }\left(  i\left(  x_{0}\right)
,r\right)  $ we do not necessarily have such representatives of $V^{2}$ which
remain in $i\left(  \overline{B}\right)  $.

With this curve $\phi$ we have $d_{\infty}\left(  \phi\left(  h\right)
,i\left(  \overline{B}\right)  \right)  =0$ and $\left(  \ref{c(h)Estim1}%
\right)  $ equals%
\[
\frac{\left\|  w+hV^{3}\left(  w\right)  -\phi\left(  h\right)  \right\|  }%
{h}=\left\|  \frac{\phi\left(  h\right)  -\phi\left(  0\right)  }{h}%
-V^{3}\left(  w\right)  \right\|  \rightarrow0
\]
as $h\rightarrow0^{+}$. Thus $\left(  \ref{VolkCondition}\right)  $ is
satisfied for $w\in i\left(  \overline{B}\right)  .$

For $w\notin i\left(  \overline{B}\right)  $ let $v\in i\left(  \overline
{B}\right)  $ be such that $d_{\infty}\left(  w,i\left(  \overline{B}\right)
\right)  =d_{\infty}\left(  w,v\right)  .$ Such a $v$ exists by the
compactness of $i\left(  \overline{B}\right)  .$ We may now apply the previous
case to $v$. Thus we have
\begin{align*}
&  \frac{d_{\infty}\left(  w+hV^{3}\left(  w\right)  ,i\left(  \overline
{B}\right)  \right)  -d_{\infty}\left(  w,i\left(  \overline{B}\right)
\right)  }{h}\\[0.1in]
&  \leq\frac{d_{\infty}\left(  w+hV^{3}\left(  w\right)  ,v+hV^{3}\left(
v\right)  \right)  +d_{\infty}\left(  v+hV^{3}\left(  v\right)  ,i\left(
\overline{B}\right)  \right)  }{h}\\[0.1in]
&  \qquad-\frac{d_{\infty}\left(  w,i\left(  \overline{B}\right)  \right)
}{h}\\[0.13in]
&  =\frac{\left\|  w-v+h\left[  V^{3}\left(  w\right)  -V^{3}\left(  v\right)
\right]  \right\|  +d_{\infty}\left(  v+hV^{3}\left(  v\right)  ,i\left(
\overline{B}\right)  \right)  }{h}\\[0.1in]
&  \qquad-\frac{d_{\infty}\left(  w,i\left(  \overline{B}\right)  \right)
}{h}\\[0.13in]
&  \leq\frac{\left\|  w-v\right\|  +hK_{1}\left\|  w-v\right\|  -d_{\infty
}\left(  w,i\left(  \overline{B}\right)  \right)  }{h}+\eta\left(  h\right)
\\[0.1in]
&  \left.  \text{(where }\eta\left(  h\right)  \rightarrow0\text{ as
}h\rightarrow0\text{)}\right. \\[0.1in]
&  =K_{1}d_{\infty}\left(  w,i\left(  \overline{B}\right)  \right)
+\eta\left(  h\right)
\end{align*}
and $\left(  \ref{VolkCondition}\right)  $ is satisfied for $w\notin i\left(
\overline{B}\right)  .$ Thus the unique solutions to $V^{3}$ with initial
conditions in $i\left(  \overline{B}\right)  $ remain in $i\left(
\overline{B}\right)  .$

Thus the solution $\sigma$ of $V^{3}$ with initial condition $i\left(
x_{0}\right)  $ exists and remains in $i\left(  \overline{B}\right)  .$ By the
continuity of $\sigma,$ there exists $\delta>0$ such that $\sigma\left(
\left[  0,\delta\right)  \right)  \subset B_{\left\|  \cdot\right\|  }\left(
i\left(  x_{0}\right)  ,r/2\right)  $ on which $V^{3}\left(  w\right)
=V^{2}\left(  w\right)  =V^{1}\left(  w\right)  $ when $w\in i\left(
\overline{B}\right)  $ so that $i^{-1}\circ\sigma:\left[  0,\delta\right)
\rightarrow X$ is a solution to $V$ with initial condition $x_{0}.$ Note that
this solution is continuous with respect to $d_{C}$ (as required), but not
necessarily with respect to $d$.
\end{proof}

\begin{lemma}
\label{fLipExten}Let $E$ be a normed vector space and let $f:E\rightarrow E$
be a $K$-Lipschitz map. For some fixed $x_{0}\in E$ and $r>0$ let $f^{\ast
}:E\rightarrow E$ to be defined as%
\[
f^{\ast}\left(  x\right)  =\left\{
\begin{array}
[c]{ll}%
f\left(  x\right)  & x\in B\left(  x_{0},r/2\right) \\
\left(  2-\frac{2}{r}\left\|  x-x_{0}\right\|  \right)  f\left(  x\right)
\quad & x\in B\left(  x_{0},r\right)  \backslash B\left(  x_{0},r/2\right) \\
0 & x\notin B\left(  x_{0},r\right)
\end{array}
\right.
\]
Then $f^{\ast}$ is Lipschitz.
\end{lemma}

\begin{proof}
Clearly $f^{\ast}$ is $K$-Lipschitz inside $B\left(  x_{0},r/2\right)  $ and
$0$-Lipschitz outside $B\left(  x_{0},r\right)  .$ Hence the analysis breaks
down into the following four cases:

Case 1. $x,y\in B\left(  x_{0},r\right)  \backslash B\left(  x_{0},r/2\right)
.$%
\begin{align*}
&  \left\|  f^{\ast}\left(  x\right)  -f^{\ast}\left(  y\right)  \right\|
\\[0.1in]
&  =\left\|  f\left(  x\right)  \left(  2-\frac{2}{r}\left\|  x-x_{0}\right\|
\right)  -f\left(  y\right)  \left(  2-\frac{2}{r}\left\|  y-x_{0}\right\|
\right)  \right\| \\[0.1in]
&  \leq2\left\|  f\left(  x\right)  -f\left(  y\right)  \right\|  +\frac{2}%
{r}\left\|  \left(  f\left(  x\right)  \left\|  x-x_{0}\right\|  -f\left(
y\right)  \left\|  y-x_{0}\right\|  \smallskip\right)  \right\| \\[0.1in]
&  =2\left\|  f\left(  x\right)  -f\left(  y\right)  \right\|  +\frac{2}%
{r}\left\|  \left(
\begin{array}
[c]{c}%
\left[  f\left(  x\right)  -f\left(  y\right)  \right]  \left\|
x-x_{0}\right\|  \quad\quad\quad\quad\quad\\
+f\left(  y\right)  \left[  \left\|  x-x_{0}\right\|  -\left\|  y-x_{0}%
\right\|  \right]
\end{array}
\right)  \right\| \\[0.1in]
&  \leq2\left\|  f\left(  x\right)  -f\left(  y\right)  \right\|  +\frac{2}%
{r}\left\|  f\left(  x\right)  -f\left(  y\right)  \right\|  \left\|
x-x_{0}\right\| \\
&  \qquad\qquad+\frac{2}{r}\left\|  f\left(  y\right)  \right\|  \left|
\left\|  x-x_{0}\right\|  -\left\|  y-x_{0}\right\|  \smallskip\right|
\\[0.1in]
&  \leq2\left\|  f\left(  x\right)  -f\left(  y\right)  \right\|  +\frac{2}%
{r}\left(  \left\|  f\left(  x\right)  -f\left(  y\right)  \right\|  \left\|
x-x_{0}\right\|  +\left\|  f\left(  y\right)  \right\|  \left\|  x-y\right\|
\smallskip\right) \\[0.1in]
&  \leq\left(  2K+\frac{2}{r}Kr+\frac{2}{r}M\right)  \left\|  x-y\right\|
\end{align*}
where $M=\sup\left\{  \left\|  f\left(  y\right)  \right\|  :y\in B\left(
x_{0},r\right)  \right\}  <\infty$ since $f$ is Lipschitz. Therefore $f^{\ast
}$ is $K_{1}$-Lipschitz with $K_{1}:=4K+\frac{2}{r}M.$

Case 2. $x\in B\left(  x_{0},r/2\right)  $ and $y\in B\left(  x_{0},r\right)
\backslash B\left(  x_{0},r/2\right)  .$

Let $z_{0}$ be a point with $\left\|  z_{0}-x_{0}\right\|  =r/2$ and
$z_{0}=t_{0}x+\left(  1-t_{0}\right)  y$ for some $0\leq t_{0}<1$ (such a
$z_{0} $ exists by continuity). Then%
\begin{align*}
\left\|  f^{\ast}\left(  x\right)  -f^{\ast}\left(  y\right)  \right\|   &
\leq\left\|  f^{\ast}\left(  x\right)  -f^{\ast}\left(  z_{0}\right)
\right\|  +\left\|  f^{\ast}\left(  z_{0}\right)  -f^{\ast}\left(  y\right)
\right\| \\
&  \leq K\left\|  x-z_{0}\right\|  +K_{1}\left\|  z_{0}-y\right\| \\
&  \leq K_{1}\left(  \left\|  x-z_{0}\right\|  +\left\|  z_{0}-y\right\|
\right) \\
&  =K_{1}\left(  \left\|  x-\left(  t_{0}x+\left(  1-t_{0}\right)  y\right)
\right\|  +\left\|  t_{0}x+\left(  1-t_{0}\right)  y-y\right\|  \right) \\
&  =K_{1}\left(  \left[  1-t_{0}\right]  \left\|  x-y\right\|  +t_{0}\left\|
x-y\right\|  \right)  =K_{1}\left\|  x-y\right\|  .
\end{align*}

Case 3. $x\in B\left(  x_{0},r\right)  \backslash B\left(  x_{0},r/2\right)  $
and $y\notin B\left(  x_{0},r\right)  .$

Similar to Case 2.

Case 4. $x\in B\left(  x_{0},r/2\right)  $ and $y\notin B\left(
x_{0},r\right)  .$

Like before, there exists a point $z_{1}$ with $\left\|  z_{1}-x_{0}\right\|
=r/2$ and $z_{1}=t_{1}x+\left(  1-t_{1}\right)  y$ for some $0<t_{1}<1.$ Then
let $z_{2}$ be a point with $\left\|  z_{2}-x_{0}\right\|  =r$ and
$z_{2}=t_{2}x+\left(  1-t_{2}\right)  y$ for some $0\leq t_{2}<1$ and
$t_{2}<t_{1}$. Then%
\begin{align*}
&  \left\|  f^{\ast}\left(  x\right)  -f^{\ast}\left(  y\right)  \right\| \\
&  \leq\left\|  f^{\ast}\left(  x\right)  -f^{\ast}\left(  z_{1}\right)
\right\|  +\left\|  f^{\ast}\left(  z_{1}\right)  -f^{\ast}\left(
z_{2}\right)  \right\|  +\left\|  f^{\ast}\left(  z_{2}\right)  -f^{\ast
}\left(  y\right)  \right\| \\
&  \leq K\left\|  x-z_{1}\right\|  +K_{1}\left\|  z_{1}-z_{2}\right\|  +0\\
&  \leq K_{1}\left(  \left\|  x-z_{1}\right\|  +\left\|  z_{1}-z_{2}\right\|
\right) \\
&  =K_{1}\left(  \left\|  x-\left(  t_{1}x+\left(  1-t_{1}\right)  y\right)
\right\|  +\left\|  t_{1}x+\left(  1-t_{1}\right)  y-\left[  t_{2}x+\left(
1-t_{2}\right)  y\right]  \right\|  \right) \\
&  =K_{1}\left(  \left[  1-t_{1}\right]  \left\|  x-y\right\|  +\left(
t_{1}-t_{2}\right)  \left\|  x-y\right\|  \right)  =K_{1}\left(
1-t_{2}\right)  \left\|  x-y\right\|
\end{align*}
so that $f^{\ast}$ is again $K_{1}$-Lipschitz.
\end{proof}

\begin{remark}
We only use the local compactness of $\left(  X,d_{C}\right)  $ at one line in
the proof of Theorem \ref{CauLips}. Perhaps a better analyst can complete the
proof assuming $\left(  X,d_{C}\right)  $ is only locally complete, in the
sense that every element $x\in X$ is contained in a complete neighborhood.
Open subsets of a complete metric space are locally complete, and it is
straightforward to show that every locally complete metric space is an open
subset of its metric completion.
\end{remark}

\begin{remark}
There should be skepticism about using $d_{C}$ instead of $d$.

For example we define the Lipschitz continuity of $V$. But in using metric
coordinates, the formulae for the vector field $V_{C}\left(  x\right)  $ will
automatically be in terms of $C.$ The most natural metric to use when checking
that the formulae are Lipschitz is $d_{C}$ and in each of the examples the
calculation is straightforward or automatic.

Still we need $X$ to be locally compact in $d_{C}$ not $d$ and our solutions
are only guaranteed to be continuous with respect to $d_{C}$. Therefore it is
important to study the connection between $\left(  X,d\right)  $ and $\left(
X,d_{C}\right)  $ which is the purpose of the following section.
\end{remark}

\section{$d_{C}$ versus $d$}

For a metric coordinate system $\left(  M,d,X,C\right)  $ the metric $d_{C}$
on $X$ defined by $\left(  \ref{dC}\right)  $can behave rather unintuitively.
E.g., there exist sequences $x_{n},y_{n}\subset\mathbb{R}^{2}$ with Euclidean
metric $d$ for which $d_{C}\left(  x_{n},y_{n}\right)  \rightarrow0$ but
$d\left(  x_{n},y_{n}\right)  \rightarrow\infty.$ As case in point, choose
$C:=\left\{  \left(  0,0\right)  ,\left(  1,0\right)  \right\}  $ and
$x_{n}:=\left(  n,2^{n}\right)  $ and $y_{n}:=\left(  -n,2^{n}\right)  .$
Therefore we are very interested in the answer to the following:

\begin{problem}
[Open question]Characterize the metric coordinate systems $\left(
M,d,X,C\right)  $ for which $\left(  X,d_{C}\right)  $ is locally compact.
\end{problem}

\begin{example}
$(X,d)$ being locally compact does not guarantee that $\left(  X,d_{C}\right)
$ is locally compact. Take $M=\mathbb{R}^{2}$ with Euclidean metric $d$,%
\[
X:=\left\{  \left(  x,y\right)  \in\mathbb{R}^{2}:y>1\text{ and }%
x\neq0\right\}  \cup\left\{  \left(  0,-1\right)  \right\}
\]
and coordinatizing set $C:=\left\{  \left(  0,0\right)  ,\left(  1,0\right)
\right\}  $. Then $\left(  X,d\right)  $ is locally compact, but $\left(
X,d_{C}\right)  $ is not locally complete (and therefore not locally compact)
since $\left(  0,-1\right)  $ has no complete neighborhood. Notice $\left(
0,-1\right)  $ is identified with the point $\left(  0,1\right)  $ by $d_{C}$.

In particular this shows there exist metric coordinate systems $\left(
M,d,X,C\right)  $ such that $\left(  X,d\right)  $ is not homeomorphic to
$\left(  X,d_{C}\right)  $.

It can be shown that $T_{\left(  0,-1\right)  }X$ is naturally identified with
a half plane (via $\pi$ from the proof of Theorem \ref{CauLips}). Thus we can
give a locally Lipschitz vector field on $\left(  X,d_{C}\right)  $ that gives
vertical translation for its flow, but no solution exists for the initial
condition $\left(  0,-1\right)  $. Thus the assumption that $X$ be locally
compact in $d_{C}$ instead of in $d$ is the correct condition for Theorem
\ref{CauLips}.

Yet another vector field $V:X\rightarrow TX$ may be given that has a diagonal
flow. Then a solution does exist with initial condition $\left(  0,-1\right)
$; but this solution is discontinuous in $\left(  X,d\right)  ,$ jumping from
$\left(  0,-1\right)  $ to $\left(  0,1\right)  .$ Such discontinuous
solutions do not exist in the case that the closed balls of $\left(
X,d\right)  $ are compact which follows from Theorem \ref{compactBalls}, below.

A similar setup shows that $\left(  X,d_{C}\right)  $ being locally compact
does not guarantee $\left(  X,d\right)  $ is locally compact. Take
$M=\mathbb{R}^{2},$%
\[
X:=\left\{  \left(  x,y\right)  \in\mathbb{R}^{2}:y\geq1\text{ and }%
x\neq0\right\}  \cup\left\{  \left(  0,1\right)  \right\}  \cup\left\{
\left(  0,y\right)  :y<-1\right\}
\]
and The coordinatizing set is $C:=\left\{  \left(  0,0\right)  ,\left(
1,0\right)  \right\}  $. $d_{C}$ ``sees'' the points $\left\{  \left(
0,y\right)  :y<-1\right\}  $ as if they were reflected across the $x$-axis.
Now $\left(  X,d_{C}\right)  $ is locally compact, but $\left(  X,d\right)  $
is not locally complete (and therefore not locally compact) since $\left(
0,1\right)  $ has no complete neighborhood.
\end{example}

Despite the pessimism of this example, we have a good beginning on answering
the open question with:

\begin{theorem}
\label{compactBalls}If all closed balls are compact in $\left(  X,d\right)  $,
then $\left(  X,d\right)  $ is homeomorphic to $\left(  X,d_{C}\right)  $, and
in particular $\left(  X,d_{C}\right)  $ is locally compact.
\end{theorem}

\begin{proof}
Since we know $d_{C}\leq d$ by $\left(  \ref{dc<d}\right)  $, we need to show
that if a sequence converges in $\left(  X,d_{C}\right)  $, then it also
converges in $\left(  X,d\right)  $. So pick a sequence $x_{n}\rightarrow x$
in $\left(  X,d_{C}\right)  $. Pick a $c\in C$. Then the sequence $d\left(
x_{n},c\right)  $ converges to $d\left(  x,c\right)  $ and is therefore
bounded, implying that the sequence $x_{n}$ is bounded in $\left(  X,d\right)
$ and thus contained in a closed, therefore compact, ball $Q.$ Now assume
$x_{n}$ does not converge in $\left(  X,d\right)  $. Then it has at least two
subsequences which converge towards different points of $Q$, say $x_{n_{i}%
}\rightarrow u$ and $x_{n_{j}}\rightarrow v$ (both with respect to $d $),
$u\neq v$. Since $d_{C}\leq d$, we know that $x_{n_{i}}\rightarrow u$ and
$x_{n_{j}}\rightarrow v$ with respect to $d_{C}$ as well. But $x_{n}$
converges to $x$ in $\left(  X,d_{C}\right)  $, so $u=x=v$ contradicting
$u\neq v$.
\end{proof}

Thus any closed subset $X$ of $\mathbb{R}^{n}$ with the Euclidean metric $d$
gives a locally compact metric space $\left(  X,d_{C}\right)  .$ This is also
true for any open subset of $\mathbb{R}^{n}$ as is proven below in Proposition
\ref{OpenDDC}.

\begin{corollary}
Let $\left(  M,d,X,C\right)  $ be a metric coordinate system and assume the
closed balls of $\left(  X,d\right)  $ are compact. Then there exist unique
solutions for any locally Lipschitz metric-coordinate vector field $V:\left(
X,d_{C}\right)  \rightarrow\left(  TX,d_{C}^{T}\right)  $. In addition all
solutions are continuous with respect to $d$.
\end{corollary}

\begin{proposition}
\label{OpenDDC}Let $\left(  M,d,X,C\right)  $ be a metric coordinate system
with $M\subset\mathbb{R}^{n}$. If $X$ is an open subset of $\mathbb{R}^{n}$
and $d$ is the Euclidean metric, then $\left(  X,d\right)  $ is homeomorphic
to $\left(  X,d_{C}\right)  .$
\end{proposition}

\begin{proof}
[Proof (Sketch)]Again pick a sequence $x_{n}\rightarrow x$ in $\left(
X,d_{C}\right)  $ and assume it does not converge in $\left(  X,d\right)  $.
This sequence is bounded in $\left(  X,d\right)  $ and therefore there exists
a subsequence $x_{n_{j}}\rightarrow y$ in $\left(  X,d\right)  $ for
$y\in\mathbb{R}^{n}$ with $y\neq x$. Also $d_{C}\left(  x,y\right)  =0$ so
$y\notin X$. Since $C$ doesn't distinguish metrically between $x$ and $y,$ $C$
must be contained in the hyperplane of $\mathbb{R}^{n}$ perpendicular to
$\overline{xy}$ through its midpoint. Further there exists $\epsilon>0$ such
that $B_{d}\left(  x,\epsilon\right)  \subset X$ since $X$ is open. Then
$B_{d}\left(  y,\epsilon\right)  $ is symmetric with respect to the hyperplane
to $B_{d}\left(  x,\epsilon\right)  $. Therefore every point in $B_{d}\left(
y,\epsilon\right)  $ has a counterpart in $B_{d}\left(  x,\epsilon\right)  $
which are not distinguished by $C$. (This uses the geometry of the Euclidean
balls.) Thus $B_{d}\left(  y,\epsilon\right)  \subset\mathbb{R}^{n}\backslash
X$, so $x_{n_{j}}\nrightarrow y$ in $\left(  X,d\right)  $ which is the
desired contradiction.
\end{proof}

The proposition relies heavily on the geometry of $\mathbb{R}^{n}$ with the
Euclidean metric as the following example shows.

\begin{example}
Consider $M:=\mathbb{R}^{2}$ with the supremum metric $d_{\infty}$. Define
$u:=\left(  1,1\right)  $, $v:=\left(  0,-1\right)  $ and
\[
X:=\left\{  x\in\mathbb{R}^{2}:d_{\infty}\left(  x,u\right)  <\tfrac{1}%
{4}\text{ or }d_{\infty}\left(  x,v\right)  <\tfrac{1}{4}\right\}
\backslash\left\{  \left(  s,t\right)  \in\mathbb{R}^{2}:s-t=1\right\}
\text{.}%
\]
$X$ is open and may be shown to be coordinatized by%
\[
C:=\left\{  \left(  0,0\right)  ,\left(  1,0\right)  ,\left(  -1,1\right)
,\left(  2,-1\right)  \right\}  .
\]
The sequence $x_{n}:=\left(  \frac{1}{n},-1-\frac{1}{n}\right)  $ for $n>4$
converges towards $u$ in $\left(  X,d_{C}\right)  $, but does not converge in
$\left(  X,d\right)  $. Therefore $\left(  X,d_{C}\right)  $ is not
homeomorphic to $\left(  X,d\right)  $.
\end{example}

\section{Invariance on metric coordinate systems}

\begin{definition}
On a metric coordinate system $\left(  M,d,X,C\right)  $ a subset $S\subset X$
is said to be \textbf{positively invariant} with respect to the
metric-coordinate vector field $V:X\rightarrow TX$ if any solution
$\sigma:\left[  0,\delta\right)  \rightarrow X$ to $V$ with initial condition
in $S$ has $\sigma\left(  t\right)  \in S$ for all $t\in\left[  0,\delta
\right)  .$
\end{definition}

We present a new version of the Nagumo Invariance Theorem on metric coordinate
systems which follows easily from work completed above.

\begin{theorem}
Let $\left(  M,d,X,C\right)  $ be a metric coordinate system for which
$\left(  X,d_{C}\right)  $ is locally compact. Let $S$ be a closed subset of
$\left(  X,d_{C}\right)  $. Let $V:\left(  X,d_{C}\right)  \rightarrow\left(
TX,d_{C}^{T}\right)  $ be a locally Lipschitz metric-coordinate vector field
tangent to $S$ in the following sense:%
\begin{gather*}
\text{For each }x\in S\text{, there exists a curve }\phi^{x}:\left[
0,\delta\right)  \rightarrow S\\
\text{which is a member of the equivalence class }V\left(  x\right)  .
\end{gather*}
Then $S$ is positively invariant with respect to $V$.
\end{theorem}

\begin{proof}
$\left(  M,d,S,C\right)  $ is a metric coordinate system and $V$ restricts to
a locally Lipschitz vector field $V|_{S}:\left(  S,d_{C}\right)
\rightarrow\left(  TS,d_{C}^{T}\right)  $ since $TS$ is naturally embedded in
$TX$. Also $S$ is locally compact, being a closed subset of a locally compact
space. Hence unique maximal solutions to $V|_{S}$ exist in $S$ by Theorem
\ref{CauLips} and coincide with solutions to $V$. If a solution $\sigma$ to
$V$ with initial condition in $S$ ever leaves $S,$ then define%
\[
t_{1}:=\sup\left\{  t:\sigma\left(  \left[  0,t\right)  \right)  \subset
S\right\}  .
\]
Since $S$ is closed in $\left(  X,d_{C}\right)  $ and $\sigma$ is continuous
with respect to $d_{C}$ we know $\sigma\left(  t_{1}\right)  \in S$. Further,
we know for the initial condition $\sigma\left(  t_{1}\right)  \in S$, the
vector field $V|_{S}$ has a solution which remains in $S$ for short time and
coincides with the solution $\sigma$ to $V$, which is a contradiction. Thus
$S$ is positively invariant.
\end{proof}

\section{Further examples and counterexamples}

With reference to Remark \ref{NonDiffCoord}, if the coordinatizing set $C$ is
a subset of $X,$ vector fields on $X$ usually cannot be nonzero and continuous
on $C$.

\begin{example}
\label{airTraffic}We work on the open half space $X=H^{3}\subset\mathbb{E}%
^{3}=M.$ For a metric coordinatizing set choose 3 points $C=\left\{
a,b,c\right\}  $ on the boundary forming a right triangle with legs of length
$d\left(  a,c\right)  =d\left(  b,c\right)  =1$ and $\overline{cb}%
\bot\overline{ca}$. E.g., if we were to use Cartesian coordinates we might
designate $H^{3} $ as the half-space with $z$-coordinate positive, and specify
$a=\left(  1,0,0\right)  $, $b=\left(  0,1,0\right)  $, and $c=\left(
0,0,0\right)  $.

For any point $x\in H^{3}$ we have $T_{x}H^{3}$ is naturally identified with
$\mathbb{R}^{3}$. Therefore it is very easy to generate vector fields in
metric coordinates. In fact a vector field $V:H^{3}\rightarrow TH^{3}$ given
by%
\[
V_{C}\left(  x\right)  :=\left(  f\left(  x_{a},x_{b},x_{c}\right)  ,g\left(
x_{a},x_{b},x_{c}\right)  ,h\left(  x_{a},x_{b},x_{c}\right)  \right)
\]
for any locally Lipschitz $f,g,h:\mathbb{R}^{3}\rightarrow\mathbb{R}$ will
always be well-defined; and by Theorem \ref{CauLips} and Proposition
\ref{OpenDDC}, $V$ is guaranteed to have unique solutions for any initial
condition in $H^{3},$ and furthermore the solutions are continuous with
respect to both metrics, $d_{C}$ and $d$. Finding the actual solutions amounts
to solving the problem as if it were a traditional vector field $V:\mathbb{R}%
^{3}\rightarrow\mathbb{R}^{3}$ and then restricting the solutions to their
domains of definition in $H^{3}.$

If we stipulate for all $w$ that $V_{b}=-V_{a}$ then all solutions will have
$\frac{d}{dt}\left[  \sigma_{a}\left(  t\right)  +\sigma_{b}\left(  t\right)
\right]  =0.$ With reference to Example \ref{H^2CoordEx}, this means the
solutions are restricted to ellipsoids with foci $a$ and $b$ since $\sigma
_{a}+\sigma_{b}$ remains constant. Alternatively, if $V_{a}=V_{b}$ then the
flows are restricted to hyperboloids with foci $a$ and $b.$ When $V_{a}=0$
then the flows are restricted to spheres with center $a$.

Define, for example, the vector field $V:H^{3}\rightarrow TH^{3}$ by%
\begin{align}
V_{a}\left(  x\right)   &  :=1\nonumber\\
V_{b}\left(  x\right)   &  :=-V_{a}\left(  x\right)  =-1\label{SphEllVF}\\
V_{c}\left(  x\right)   &  :=0.\nonumber
\end{align}
For an initial condition $x$ in metric coordinates $x_{C}=\left(  x_{a}%
,x_{b},x_{c}\right)  \in H^{3},$ the solution $\sigma$ follows the
intersection of the ellipsoid with foci $a$ and $b$ which touches $x$ and the
sphere centered at $c$ which touches $x$. The formula is found by regular
integration of $\left(  \ref{SphEllVF}\right)  $ to be%
\[
\sigma_{C}\left(  t\right)  =\left(  x_{a}+t,x_{b}-t,x_{c}\right)
\]
in metric coordinates. One particular solution $\sigma$ is graphed in Figure
\ref{SolutionFigure1}. That the graph of $\sigma$ is given by the intersection
of an ellipsoid and sphere is illustrated in Figure \ref{IntersectionFigure2}.%
%TCIMACRO{\FRAME{ftbphFU}{2.834in}{2.1326in}{0pt}{\Qcb{Figure 1: Solution
%curve}}{\Qlb{SolutionFigure1}}{fig1.png}%
%width 2.834in;  height 2.1326in;  depth 0pt;  original-width 2.7916in;
%original-height 2.0937in;  cropleft "0";  croptop "1";  cropright "1";
%cropbottom "0";  filename 'Fig1.png';file-properties "XNPEU";}}}%
%BeginExpansion
%%\begin{figure}
%%[ptbh]
%%\begin{center}
%%\includegraphics[
%%natheight=2.093700in,
%%natwidth=2.791600in,
%%height=2.1326in,
%%width=2.834in
%%]%
%%{Fig1.eps}%
%%\caption{Figure 1: Solution curve}%
%%\label{SolutionFigure1}%
%%\end{center}
%%\end{figure}
%EndExpansion%
%TCIMACRO{\FRAME{ftbphFU}{3.6409in}{2.866in}{0pt}{\Qcb{Figure 2: Intersection
%of sphere and ellipsoid}}{\Qlb{IntersectionFigure2}}{fig2.png}%
%width 3.6409in;  height 2.866in;  depth 0pt;  original-width 3.5942in;
%original-height 2.8228in;  cropleft "0";  croptop "1";  cropright "1";
%cropbottom "0";  filename 'Fig2.png';file-properties "XNPEU";}}}%
%BeginExpansion
%%\begin{figure}
%%[ptbhptbh]
%%\begin{center}
%%\includegraphics[
%%natheight=2.822800in,
%%natwidth=3.594200in,
%%height=2.866in,
%%width=3.6409in
%%]%
%%{Fig2.eps}%
%%\caption{Figure 2: Intersection of sphere and ellipsoid}%
%%\label{IntersectionFigure2}%
%%\end{center}
%%\end{figure}
%EndExpansion

\begin{center}  \label{SolutionFigure1}
\epsfig{file=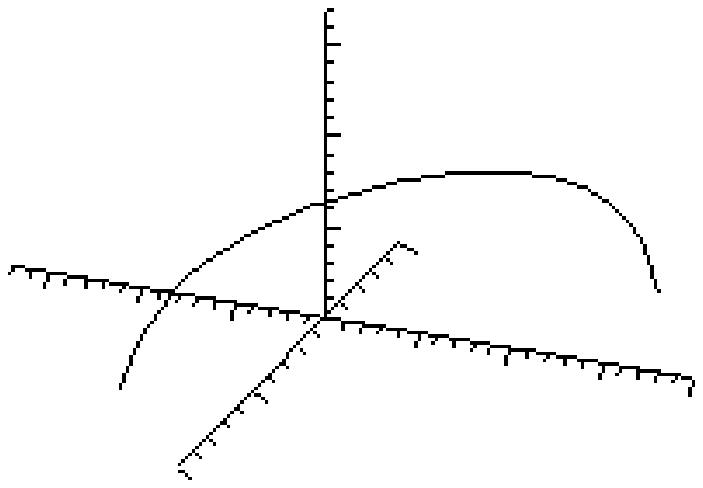,height=2in,width=2.8in} \\ {\sc Figure 1:
Solution Curve}
\end{center}

\begin{center}  \label{IntersectionFigure2}
\epsfig{file=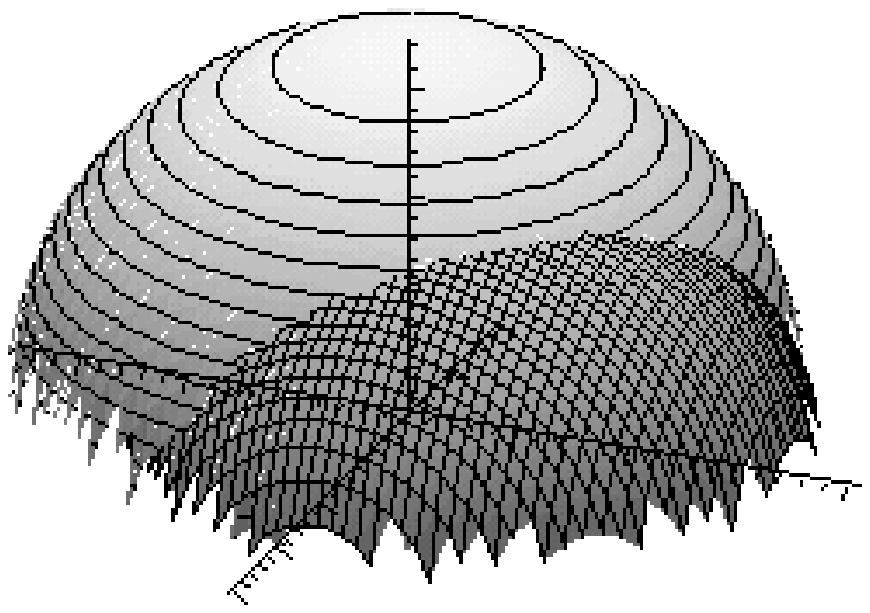,height=2.9in,width=3.7in} \\ {\sc Figure 2:
Intersection of sphere and ellipsoid}
\end{center}

Alternatively examine the metric-coordinate vector field given by%
\begin{align*}
V_{a}\left(  x\right)   &  :=1\\
V_{b}\left(  x\right)   &  :=-V_{a}\left(  x\right)  =-1\\
V_{a}\left(  x\right)   &  :=V_{b}\left(  x\right)  =1.
\end{align*}
Then for an initial condition $x\in H^{3},$ the solution $\sigma$ follows the
intersection of the ellipsoid with foci $a$ and $b$ which touches $x$ and the
hyperboloid with foci $a$ and $b$ touching $x.$ The formula is simply%
\[
\sigma_{C}\left(  t\right)  =\left(  x_{a}+t,x_{b}-t,x_{c}-t\right)  .
\]

On the boundary metric-coordinate vector fields are not so easily generated
since the tangent spaces are not all of $\mathbb{R}^{3}.$
\end{example}

\begin{example}
\label{ExS^2}Let $M=X:=S^{2}$ be the Euclidean sphere with radius 1 and
intrinsic metric $d\left(  x,y\right)  $ given by the length of a shortest
geodesic connecting $x$ and $y.$ Metrically coordinatize $S^{2}$ with
$C:=\left\{  a,b,c\right\}  $ where the three points are chosen so that
\[
d\left(  a,b\right)  =d\left(  b,c\right)  =d\left(  c,a\right)  =\frac{\pi
}{2}.
\]
We wish to have solutions follow hyperbolic paths on $S^{2}$ with foci $a$ and
$b.$ We thus need to define $V:S^{2}\rightarrow TS^{2}\subset\mathbb{R}^{3}$
(by suppressing the notation of $\pi$) with%
\begin{align*}
V_{a}\left(  x\right)   &  :=f\left(  d\left(  x,a\right)  ,d\left(
x,b\right)  ,d\left(  x,c\right)  \right)  =f\left(  x_{C}\right) \\
V_{b}\left(  x\right)   &  :=g\left(  d\left(  x,a\right)  ,d\left(
x,b\right)  ,d\left(  x,c\right)  \right)  =g\left(  x_{C}\right) \\
V_{c}\left(  x\right)   &  :=h\left(  d\left(  x,a\right)  ,d\left(
x,b\right)  ,d\left(  x,c\right)  \right)  =h\left(  x_{C}\right)
\end{align*}
where the functions $f,g,h:\mathbb{R}^{3}\mathbb{\rightarrow R}$ are Lipschitz
with $f=g$. Some further conditions on $f,g,$ and $h$ are necessary to get a
\textit{bona-fide} map into $TS^{2}.$ The hyperbolic paths of the solutions to
$V$ will be perpendicular to the great circle $S^{1}$ through $a$ and $b$
which is given by%
\[
S^{1}:=\left\{  x:d(x,c)=\frac{\pi}{2}\right\}
\]
Thus the rate of change of the distance from $a$ to a solution curve $\sigma$
will be $0$ if $\sigma$ passes through $S^{1};$ i.e., for $\sigma\left(
t\right)  =x\in S^{1}$ we need $\sigma_{a}^{+}\left(  t\right)  =V_{a}\left(
x\right)  =0$. Thus make%
\begin{equation}
f\left(  u,v,w\right)  =0\quad\text{ for \quad}w=d\left(  x,c\right)
=\frac{\pi}{2}. \label{fCond1}%
\end{equation}
E.g.,
\[
f\left(  u,v,w\right)  :=\left(  \frac{\pi}{2}-w\right)
\]
with $g=f$, then $h$ is determined by $f$ and $g$ and by continuity. This
gives a hyperbolic flow with foci $a$ and $b.$ The direction of the flow is
determined by
\begin{align}
f\left(  u,v,w\right)   &  >0\quad\text{for}\quad w>\frac{\pi}{2}%
\quad\text{and}\label{fCond2}\\
f\left(  u,v,w\right)   &  <0\quad\text{for}\quad w<\frac{\pi}{2}\text{.}
\label{fCond3}%
\end{align}
Thus on the hemisphere bounded by $S^{1}$ containing $c$ the flow is toward
$a$ and $b.$ On the complimentary hemisphere solutions move away from $a$ and $b.$
\end{example}

\begin{example}
The fact that the metrics $d_{C}$ and $d$ are not always equivalent may be
exploited to give discontinuous flows. Let $M:=\mathbb{R}^{2}$ with Euclidean
metric $d$ and let $X$ consist of the infinite strips%
\[
X:=\mathbb{R}\times\left(  \left(  0,1\right)  \underset{n\in\mathbb{N}}{\cup
}\left(  \left[  2n,2n+1\right)  \cup\left(  -2n-2,-2n-1\right]  \right)
\right)  .
\]
Let $C:=\left\{  a=\left(  0,0\right)  ,b=\left(  1,0\right)  \right\}  $, so
that $\left(  X,d_{C}\right)  $ may be identified with an open half plane and
is locally compact. A locally Lipschitz vector field $V:\left(  X,d_{C}%
\right)  \rightarrow\left(  TX,d_{C}^{T}\right)  $ is given by $V_{a}\left(
x\right)  :=1$ and $V_{b}\left(  x\right)  :=1$ so that as before we get a
hyperbolic flow on the half plane. Such solutions, however, are discontinuous
with respect to $d$, jumping from strip to strip.
\end{example}

\begin{example}
[Observer dependence of tangent spaces]On general metric spaces the tangent
space at a point fundamentally depends on the choice of metric coordinatizing
set. The metric space from Example \ref{CoordDep} is an obvious candidate to
consider and does give the result we seek: with respect to $\left\{  a=\left(
-2,0\right)  \right\}  $ the tangent space at $\left(  0,0\right)  $ may be
identified with $\mathbb{R}$ (via the weak contraction $\pi$ where $\pi\left(
\left[  \phi\right]  \right)  :=\phi_{C}^{+}\left(  0\right)  $) and with
respect to $\left\{  b=\left(  1,1\right)  \right\}  $ the tangent space
$T_{\left(  0,0\right)  }X$ is $\mathbb{R}^{-}:=\left(  -\infty,0\right]  $.
The discrepancy arises because a curve issuing from $\left(  0,0\right)  $ in
the direction of the point $\left(  -1,1\right)  $ with finite Euclidean speed
will be tangent to the circle with center $b=\left(  1,1\right)  $ and radius
$\sqrt{2}.$ Hence, the rate of change of distance will be zero with respect to
the metric coordinate $b=\left(  1,1\right)  ;$ therefore there is no positive
representative in $T_{\left(  0,0\right)  }X$. With respect to $a=\left(
-2,0\right)  $ however, a curve can issue from $\left(  0,0\right)  $ with
positive or negative $y$-metric-coordinate derivative.

If we were to allow curves with infinite speed at $t=0$ to represent members
of $T_{\left(  0,0\right)  }X$ we could recover all of $\mathbb{R}$ for the
tangent space with respect to $\left\{  \left(  1,1\right)  \right\}  $. For
example, the curve $\phi\left(  t\right)  :=\left(  -\sqrt{t},\sqrt{t}\right)
$ has
\[
\phi_{b}\left(  t\right)  =d\left(  \phi\left(  t\right)  ,\left(  1,1\right)
\right)  =\left(  \sqrt{\left(  -\sqrt{t}-1\right)  ^{2}+\left(  \sqrt
{t}-1\right)  ^{2}}\right)  =\sqrt{2t+2}%
\]
so that $\phi_{b}^{\prime}\left(  0\right)  =\frac{1}{\sqrt{2}}$. Then the
tangent spaces in this example would be topologically equivalent with respect
to different coordinate systems. Still, there exist metric coordinate systems
where this fails to patch up the disparity as in the following:
\end{example}

\begin{example}
In $M:=\mathbb{R}^{2}$ with the Euclidean metric $d$, choose%
\[
X:=\left\{  \left(  t,t\sin\left(  \frac{1}{t}\right)  \right)  \in
\mathbb{R}^{2}:0<\left|  t\right|  \leq1\right\}  \cup\left\{  \left(
0,0\right)  \right\}  \text{.}%
\]
Let $c_{0}=\left(  0,1\right)  ,$ $c_{1}=\left(  1,1\right)  ,$ $c_{2}=\left(
-1,0\right)  ,$ $c_{3}=\left(  1,0\right)  $, and $x:=$ $\left(  0,0\right)
.$ Notice that $C_{1}:=\left\{  c_{0},c_{1}\right\}  $ and $C_{2}:=\left\{
c_{2},c_{3}\right\}  $ each metrically coordinatize $X$.

Notationally use $T_{x}^{C_{i}}X$ to denote the tangent space of $X$ at $x$
relative to $C_{i}.$ Without providing the voluminous details, we claim that
$T_{x}^{C_{1}}X$ is the singleton $0$ while $T_{x}^{C_{2}}X$ is $\mathbb{R}.$
Infinite or even $0$ speed reparametrizations will not recover any other
elements of $T_{x}^{C_{1}}X.$
\end{example}

\begin{example}
To amplify the last example, we show that there are dynamics describable with
a vector field with respect to one coordinatizing set which are not achievable
by any vector field with respect to another coordinatizing set.

Consider $M:=\mathbb{R}^{2}$ with the Euclidean metric $d$ and $X:=\left\{
\left(  x,\left|  x\right|  \right)  :-1\leq x\leq1\right\}  .$ Then
$C_{1}:=\left\{  \left(  -1,1\right)  \right\}  $ and $C_{2}:=\left\{  \left(
1,1\right)  \right\}  $ each coordinatize $X$. The curve%
\[
\phi\left(  t\right)  :=\left\{
\begin{array}
[c]{ll}%
\left(  2\sqrt{t},2\sqrt{t}\right)  & \quad t\geq0\\
\left(  t,\left|  t\right|  \right)  & \quad t<0
\end{array}
\right.
\]
has bounded derivative with respect to $C_{1}$ but the metric coordinate
derivative of $\phi$ with respect to $C_{2}$ does not exist at $t=0.$ The
derivative of $\phi$ then yields a vector field with respect to $C_{1}$ giving
dynamics which cannot be described with respect to $C_{2}.$
\end{example}

\section{Open questions and future directions}

Is there a canonical method for coordinatizing a metric space with a minimum
number of metric coordinates? A simpler question is: does every metric space
$\left(  X,d\right)  $ have a discrete metric coordinatizing subset $C$? Minor
headway on this latter question is given by:

\begin{remark}
If $\left(  M,d,X,C\right)  $ is a metric coordinate system and $c$ is an
accumulation point of $C$, then $C\backslash\left\{  c\right\}  $ is still a
metric coordinatizing set for $X$.
\end{remark}

\begin{proof}
If $c_{i}\rightarrow c$ and $d\left(  x,c_{i}\right)  =d\left(  y,c_{i}%
\right)  $ for all $i$ then $d\left(  x,c\right)  =d\left(  y,c\right)  $ by
continuity of the metric.
\end{proof}

\begin{example}
By the above remark, we can remove accumulation points, one at a time, from
any metric coordinatizing set. However, it may not be possible to remove all
of them. More succinctly, not every coordinatizing set has a discrete
coordinatizing subset. E.g., the closed upper half plane in $\mathbb{R}^{2}$
given by $\mathbb{R\times R}^{+}=\left\{  \left(  x,y\right)  :y\geq0\right\}
$ with the metric $d_{\infty}$ from Example \ref{InfinityMetric} is
coordinatized by the line $C:=\left\{  \left(  x,0\right)  :x\in
\mathbb{R}\right\}  .$ But any discrete subset of $C$ fails to coordinatize
$\mathbb{R\times R}^{+}.$ In fact any coordinatizing subset must be dense in
$C$. The open question remains, however, since there does exist a discrete
coordinatizing set on $\left(  \mathbb{R\times R}^{+},d_{\infty}\right)  ,$
namely the one from Example \ref{InfinityMetric}.
\end{example}

We also have the open question from the end of Section \ref{CauLipSection}:
characterize the metric coordinate systems $\left(  M,d,X,C\right)  $ for
which $\left(  X,d_{C}\right)  $ is locally compact. A version of Theorem
\ref{CauLips} which does not require local compactness is highly desirable.
The imagined condition is that $\left(  X,d_{C}\right)  $ is locally complete.
Thus we would also be pleased with a characterization of the metric coordinate
systems $\left(  M,d,X,C\right)  $ for which $\left(  X,d_{C}\right)  $ is
locally complete.

Next, in a metric coordinate system a new measure of the dimension of a metric
space presents itself.

\begin{definition}
Let $I$ be a cardinal number. The metric space $\left(  M,d\right)  $ is
\textbf{locally }$I$\textbf{-coordinatizable} if for each $x\in M$ there
exists a neighborhood $X$ and a set $C\subset M$ of cardinality $I$ which
metrically coordinatizes $X.$The smallest such cardinal number $I$ is called
the \textbf{metric-coordinate dimension}.
\end{definition}

Metric-coordinate dimension is not a homeomorphic invariant. For example the
Koch curve is homeomorphic to $\mathbb{R}$ but has metric-coordinate dimension 2.

\begin{conjecture}
Metric-coordinate dimension is a lipeomorphic invariant\footnote{A
lipeomorphism between metric spaces $M$ and $N$ is a bijective Lipschitz map
between $M$ and $N$ with Lipschitz inverse.}.
\end{conjecture}

Next, what is the most appropriate definition for metric-coordinate-wise
differentiability of maps between metric spaces? Which brings us to question
what conditions give an Inverse Function Theorem on coordinatized metric
spaces (this has been done before on metric spaces using the structure of
``mutations'', \cite{Aubin}).

Higher order derivatives are obviously defined with $\phi_{C}^{\prime\prime
}\left(  t\right)  :=\left(  \frac{d^{2}}{dt^{2}}\phi_{c}\left(  t\right)
\right)  $. How do we analyze higher-order differential equations?

The directional derivative $D_{V}f$ of a function $f:X\rightarrow\mathbb{R}$
on a metric coordinate system in the direction of a metric-coordinate vector
field $V$ can be defined as%
\begin{equation}
D_{V}f\left(  x\right)  :=\underset{h\rightarrow0^{+}}{\lim}\frac{f\left(
\phi\left(  h\right)  \right)  -f\left(  x\right)  }{h} \label{DVWDef}%
\end{equation}
assuming the limit exists and does not depend on the representative
$\phi:\left[  0,\delta\right)  \rightarrow X$ of the equivalence class
$V\left(  x\right)  $. This notion is useful in analyzing the qualitative
dynamics of metric-coordinate vector fields using Lyapunov functions $f$ as
will be demonstrated in a forthcoming paper. Perhaps we can also use the
directional derivative to analyze extrema and extract the Lagrange multiplier
method for constraints. Certainly the fundamental theorem of line integrals
should have an expression on metric spaces with $D_{V}f$. What can be made of
Stokes' Theorem?

PDE's on metric spaces should be possible to formulate with these directional derivatives.

When we consider non-autonomous (i.e., time-dependent) metric-coordinate
vector fields, we might allow the location of the metric coordinatizing points
$c$ to change in time as well; giving us bonus descriptive power not available
with Cartesian coordinates.

Finally one might abandon the goal of finding coordinatizing sets. Begin with
any set $C\subset M$ and define the quotient space $X/\sim$ with the
equivalence relation $x\sim y$ if $d\left(  x,c\right)  =d\left(  y,c\right)
$ for all $c\in C.$ Then work in the metric space $\left(  X/\sim
,d_{C}\right)  $ which is identified with a subset of $\mathbb{R}_{b}^{C}$.

\end{document}